\documentclass[12pt]{article}
\usepackage{amsfonts}
\usepackage{amssymb,amsmath,amsthm,hyperref}

\usepackage{graphics,graphpap}
\usepackage{supertabular}

\theoremstyle{plain}

\def\mg{\mathfrak{g}}
\def\mh{\mathfrak{h}}

\def\mb{\mathfrak{b}}

\def\mn{\mathfrak{n}}
\def\mc{\mathfrak{c}}
\def\tb{\tilde{\mathfrak{b}}}
\def\tn{\tilde{\mathfrak{n}}}

\newtheorem{thm}{Theorem}[section]
\newtheorem{lem}[thm]{Lemma}

\theoremstyle{definition}

\numberwithin{equation}{section}

\def\a{\alpha}
\def\b{\beta}
\def\d{\delta}
\def\D{\Delta}

\def\si{\sigma}

\def\mb{\mathfrak{b}}
\def\mc{\mathfrak{c}}
\def\mg{\mathfrak{g}}
\def\mh{\mathfrak{h}}

\def\ad{{\text{ad}}}

\def\dim{{\text{dim}}}
\def\e{\varepsilon}

\def\N{\mathcal{N}}

\def\Z{\mathcal{Z}}
\def\C{\mathcal{C}}

\def\leftitem#1{\item{\hbox to\parindent{\enspace#1\hfill}}}

\def\dim{\operatorname{dim}}

\def\deg{\operatorname{deg}}
\def\min{\operatorname{min}}

\def\Z{\Bbb Z}\def\N{\Bbb N}

  \def\D{{\it\Delta}} \def\b{\beta}
  \def\a{\alpha}
 \def\d{\delta}  
  \def\D{\Delta} 
\def\si{\sigma}

\def\C{\Bbb C}

\def\tilde{\widetilde}

\def\Der{{\rm{Der}}}
\def\Aut{{\rm{Aut}}}

\begin{document}

\centerline {{\Large\bf Automorphisms and derivations of Borel subalgebras
}}

\centerline {{\Large\bf and their nilradicals in Kac-Moody algebras}
}

\vskip 15pt
\par
\centerline{
Jun Morita\footnote{Partially
supported by Monkasho Kakenhi of Japan}
and Kaiming Zhao\footnote{Partially
supported by NSERC of Canada and NSF of China}
}
\par
\vskip 5pt

\begin{abstract} In this paper, we determine derivations of
Borel subalgebras and their derived subalgebras called
nilradicals, in Kac-Moody algebras (and contragredient Lie
algebras) over any field of characteristic $0$; and we also
determine automorphisms of those subalgebras in symmetrizable
Kac-Moody algebras. The results solve a conjecture posed by
R.~V.~Moody about 30 years ago which generalizes a result by
B.~Kostant and which is discussed by A.~Fialowski using Lie
algebra cohomology in case of affine type.

\vspace{5mm}

\noindent\emph{Keywords}: Kac-Moody algebra, Borel subalgebra, nilradical, derivation,
automorphism

\noindent\emph{Mathematics Subject Classification} 2000: 17B05,
 17B30, 17B65, 17B67.

\indent \hskip .3cm
\end{abstract}

\section{\textbf{Introduction}}

A {\it generalized Cartan matrix} (GCM) is a matrix of integers $A =
(a_{ij})_{i,j=0}^l$ satisfying $a_{ii} = 2$ for all $i$;
$a_{ij}\le0$ if $i \ne j$; and $a_{ji}=0$ iff ~$a_{ij}=0$. $A$ is
called {\it symmetrizable} if $DA$ is symmetric for some
nondegenerate rational diagonal matrix $D$. We fix   a generalized
Cartan matrix $A$ in this paper, assumed for simplicity to be
indecomposable. Fix a base field $F$ of characteristic zero. Let
$(\mh,\Pi,\Pi^\vee)$ be a realization of $A$, where $\Pi = \{
\a_0,\a_1,\cdots, \a_l\} $ and $\Pi^\vee = \{
\a_0^\vee,\a_1^\vee,\cdots, \a_l^\vee\}$ with $ \a_i(\a_j^\vee) =
a_{ji}$ (cf.~\cite{Kac},\cite{MP}).

The {\it contragredient Lie algebra} $L(A)$ is the Lie algebra over
$F$ generated by $\mh$ and symbols $e_i$ and $f_i$ ($i \in I=\{
0,1,2,\cdots,l \}$) with defining relations:\ $[\mh,\mh] = 0;\
[e_i,f_j] = \delta_{i,j}\a_i^\vee\ (i,j \in I);\ [h,e_i] =
\a_i(h)e_i,\ [h,f_i] =-\a_i (h)f_i\ (h\in \mh, i \in I) ;\ (\ad
e_i)^{1-a_{ij}} (e_j) = 0 = (\ad f_i)^{1-a_{ij}} (f_j)$ for $i\ne j\
(i,j \in I)$. Let us denote  by $\tn^+$ the subalgebra generated by
$e_i, i\in I$, and $\tb^+=\mh\oplus\tn^+$.

The {\it Kac-Moody algebra} $\mg = \mg(A)$ is the quotient Lie
algebra $L(A)/J(A)$ where $J(A)$ is the maximal ideal of $L(A)$
trivially intersecting with $\mh$. When $A$ is symmetrizable, $J(A)$
is actually $0$ (cf.~\cite{GK},\cite{MP}). We have the canonical
embedding $\mh\subset \mg$  and linearly independent Chevalley
generators $e_i,f_i$ $(i \in I)$ for the derived algebra $\mg'$ of
$\mg$. The center $\mc$ of $\mg$ lies in $\mh': = \mh\cap \mg' =
\sum_{i=1}^l F\a_i^\vee$. Every ideal of $\mg$ contains $\mg'$ or is
contained in $\mc$. Define an involution $\omega$ of $\mg$ by
requiring $\omega(e_i) = -f_i,\, \omega(f_i) = -e_i$ and $\omega(h)
= -h$ for all $i \in I$ and $h \in \mh$. Let $\mn^+$ be the
subalgebra of $\mg$ generated by the $e_i$ $(i \in I)$, and put
$\mb^+=\mh\oplus\mn^+$, called the standard Borel subalgebra of
$\mg$. Set $\mb^-=\omega (\mb^+)$ and $\mn^- = \omega(\mn^+)$. Then,
we have $\mg = \mn^-\oplus\mh\oplus\mn^+$ (a triangular
decomposition). The subalgebras $\mn^\pm$ are sometimes called the
nilradicals of $\mb^\pm$ respectively. If $A$ is symmetrizable, we
know that $\mb^+=\tb^+$ , $\mn^+=\tn^+$.

 In this paper, we will
determine derivations $\Der(\mb^+)$, $\Der(\tb^+)$ , $\Der(\mn^+)$
and $\Der(\tn^+)$ in general case, and we will also determine
automorphisms $\Aut(\mb^+)$ and $\Aut(\mn^+)$ in symmetrizable
case. Derivations of Lie algebras can be discussed using Lie
algebra cohomology (or sometimes homology) theory
(cf.~\cite{FF},\cite{HS},\cite{Ko}). We do not use cohomology
method to compute $\Der(\mb^+)$, $\Der(\tb^+)$ , $\Der(\mn^+)$  or
$\Der(\tn^+)$ in this paper since some exact sequences do not seem
to carry enough information for indefinite type Kac-Moody
algebras. Derivations $\Der(\mb^+)$ and $\Der(\mn^+)$ in case of
finite type are already known (cf.~\cite{LL1},\cite{LL2}). Also
derivations $\Der(\mn^+)$ in case of affine type are already
obtained (cf.~\cite{F1},\cite{F2}). Automorphisms $\Aut(\mb^+)$
and $\Aut(\mn^+)$ might be known in case of finite type, but we
could not find any reference. About 30 years ago, R.~V.~Moody gave
a question about the structure of $\Der(\tn^+)$. As is obtained in
the above references, there exist derivations in $\Der(\tn^+)$
other than $\ad(\tb^+)|_{\tn^+}$ for finite type and for affine
type. Therefore the remaining problem is to study:\par
\vspace{-0.15cm} \ \textquoteleft \textquoteleft {\it It is
conceivable that} $\Der(\tn^+)$ {\it is equal to}
$\ad(\tb^+)|_{\tn^+}$ {\it when $A$ is of indefinite type},''\par
\vspace{-0.15cm} \noindent which we call {\it Moody's conjecture}
(cf.~\cite{M}). Here we will give a complete answer, which means
that Moody's conjecture is true, i.e. $\Der(\tn^+) =
\ad(\tb^+)|_{\tn^+}$ (and $\Der(\mn^+) = \ad(\mb^+)|_{\mn^+}$) for
indefinite case. Note that $\Der(\mg)$ is known
(cf.~\cite{B1},\cite{B2},\cite{Sa}, also see Theorem 5.2) in
general, and note that $\Aut(\mg)$ is also known in symmetrizable
case (cf.~\cite{MP},\cite{PK}, also see Theorem 2.2). As an
advanced result on indefinite Kac-Moody algebras, we refer to
\cite{BKM} for example.

We will denote the sets of integers, positive integers,
non-negative integers, rational numbers and complex numbers by
$\Z$, $\N$, $\Z^+$, $\Bbb Q$ and $\Bbb C$ respectively.

We have $\mg =\oplus_{\a\in\mh^*}\, \mg_{\a}$, where $\mg_{\a}=
\{x \in\mg \,|\,[h,x] = \a(h)x \ (\forall \,h \in\mh)\}$. Put $Q =
\oplus_{i \in I}\, \Z\a_i$ and $Q^+ = \sum_{i \in I}\, \Z^+\a_i$,
and define a partial order on $\mh^*$ by: $\a\le\b$ if $\b-\a\in
Q^+$. A root (resp. positive root) is an element of $\Delta = \{\a
\in\mh^*\,|\,\a \ne 0,\, \mg_\a \ne 0\}$ (resp. $\Delta_+: =
\Delta\cap Q^+$). We have $\mh=\mg_0$, $\mn^{\pm} =
\oplus_{\a\in\Delta_+}\mg_{\pm\a}$ and $\mg = \mh \oplus
(\oplus_{\a \in \D}\, \mg_\a)$, called the {\it root space
decomposition} with respect to $\mh$.

Let $s_i$ be an invertible endomorphism of $\mh^*$ defined by
$s_i(\mu) = \mu - \mu(\a_i^\vee)\a_i$ for all $\mu \in \mh^*$. Then,
the {\it Weyl group} $W$ of $\mg$ is defined to be the subgroup of
$GL(\mh^*)$ generated by $s_i$ for all $i \in I$.

Let $\Aut(A)$ be the group of all permutations $\si$ on $I$
satisfying $a_{\si(i)\si(j)} = a_{ij}$. Then, $\Aut(A)$ is called
the {\it Dynkin diagram automorphism} group of $A$, which is
regarded as a subgroup of $\Aut(\mg')$ (or $\Aut(L)$) by requiring
$\si(e_i) = e_{\si(i)}$, $\si(f_i) = f_{\si(i)}$ for $\si \in
\Aut(A)$.

If $A$ is symmetrizable, we choose a nondegenerate $\mg$-invariant
symmetric $F$-bilinear form $(.|.)$ on $\mg$ such that $(\a_i
^\vee|\a_i^\vee)$ is positive in $\mathbb{Q} \subset F$ for all $i
\in I$. Then $(.|.)$  is nondegenerate on $\mh$, and, hence,
induces a $W$-{\it invariant form} $(.|.)$ on $\mh^*$. One has
$(\a|\a)>0$ in $\mathbb{Q} \subset F$ if $\a\in\Delta^{re}$. We
also note that $(.|.)$  induces a nondegenerate form on
$\mg'/\mc$, which is unique up to multiples.

We now recall the adjoint {\it Kac-Moody group} $G$ associated to
$\mg$ (see \cite{Ku},\cite{MP},\cite{PK},\cite{Ti}). For each real root $\a \in
\D^{\rm re} = W \Pi = \{ w(\a)\,|\, w \in W,\, \a \in \Pi \}$ and
a root vector $x_\a \in \mg_\a$, we see that $\exp(\ad x_\a)$
gives an automorphism of $\mg$. Let $U_\a$ be the subgroup of the
automorphism group $\Aut(\mg)$ of $\mg$ generated by $\exp (\ad
x_\a)$ for all $x_\a \in \mg_\a$. Each $\xi \in T = {\rm
hom}_\Bbb{Z}(Q,F^\times)$ induces an automorphism of $\mg$ by
$\xi(x) = \xi(\a)x$ for all $x \in \mg_\a$ and $\a \in \D \cup \{
0 \}$. Then, the adjoint Kac-Moody group $G$ is defined to be the
subgroup of $\Aut(\mg)$ generated by $T$ and $U_\a$ for all $\a
\in \D^{\rm re}$. Note that $T$ normalizes $U_\a$. Put $\D^{\rm
re}_+ = \D^{\rm re} \cap \D_+$. Let $U_+$ (resp. $U_-$) be the
subgroup of $G$ generated by $U_\a$ (resp. $U_{-\a}$) for all $\a
\in \D^{\rm re}_+$, and put $B_\pm = TU_\pm$. Set $N = N_G(T)$,
the normalizer of $T$ in $G$. Then, we have Tits systems
$(G,B_\pm,N,S)$ and we see $N/T \simeq W$. Furthermore we have
Bruhat decompositions $G = U_\pm N U_\pm$ and Birkhoff
decompositions $G = U_\mp N U_\pm$. Put $\D^{\rm im} = \D
\setminus \D^{\rm re}$ (imaginary roots). When $F=\C$,
we see $T=\exp(\ad\ \mh)$.

Since both adjoint representations of $L(A)$ and $\mg(A)$ are
integrable,  $L(A)$ and $\mg(A)$ share a lot of properties. For
example they have the same root system and the isomorphic Weyl
group. The only difference between $L(A)$ and $\mg(A)$ is the root
spaces corresponding to imaginary roots.

 We assume that $A$ is symmetrizable
in Sections 2 and 4. In Section 2, we use Peterson-Kac's
Conjugation Theorem to determine the automorphism group
$\Aut(\mb^+)$ of $\mb^+$ (see Theorem 2.1). In particular, if $A$
is nonsingular, then $\Aut(\mb^+) = \Aut(A) \ltimes B_+$. In
Sections 3 and 5, $A$ is not necessarily symmetrizable. In Section
3, we establish an elementary method and some techniques to
determine the derivation algebras $\Der(\mn^+)$ and $\Der(\tn^+)$
for all cases (finite, affine and indefinite types) which gives a
complete answer to Moody's conjecture, and which recovers the
result by Kostant (cf.~\cite{Ko}) for finite type and the result
by A.~Fialowski (cf.~\cite{F2}) for affine case (see Theorem 3.4).
The main tricks are the technical computations using combinatorics
and calculus formulas to obtain (3.4), (3.7), (3.15) and (3.16).
In Section 4, we explicitly determine the automorphism group
$\Aut(\mn^+)$ of $\mn^+$, using the results established in Section
3 (see Theorems 4.1 and 4.2). This problem was discussed by Moody
(cf.~\cite{M}) with partial results. In Section 5, we determine
the derivation algebras $\Der(\mb^+)$ and $\Der(\tb^+)$ (see
Theorem 5.1). We reach $\Der(\mb^+) = \hom(\mh,\mc) \oplus
\ad(\mb^+)$ and $\Der(\tb^+) = \hom(\mh,\mc) \oplus \ad(\tb^+)$.
From our results and proofs, we again reveal the close
relationship between the derivation algebra and the automorphism
group of a Lie algebra (cf.~Sect.5, Remark 6).

 It is striking that the derivation algebras
$\Der(\mn^+)$ for finite type, affine type and indefinite type
Kac-Moody algebras are totally different, but the automorphism
groups $\Aut(\mn^+)$ are only different between finte type and
infinite type Kac-Moody algebras, while $\Der(\mb^+)$ and
$\Aut(\mb^+)$ have the same formulas for all (symmetrizable)
Kac-Moody algebras.

\section{\textbf{ Automorphisms of $\mb^+$}}

In this section, a generalized Cartan matrix $A$ is assumed to be
symmetrizable. Let $\dim \mh = 2 (l +1)- m = (l+1) + m'$, where
$m$ is the matrix rank of $A$ and $m' = (l+1) - m$ is the
nullity (or matrix corank) of $A$. We put $R = \oplus_{i \in I}
F\a_i$, and choose its complement $S$ in $\mh^*$, that is, $\mh^*
= R \oplus S$. We fix a basis $\{ v_1,\ldots,v_{m'} \}$ of $S$.
Then, we have the corresponding dual basis $\{
\a_0^*,\a_1^*,\ldots,\a_l^*,v_1^*,\ldots, v_{m'}^* \}$ of $\mh =
\mh^{**}$, i.e., $\a_i(\a_j^*)=\delta_{i,j}$ and
$v_i(v_j^*)=\delta_{i,j}$, further, $ \a_i(v_j^*)=v_i(\a_j^*)=0$.
Then $\mc = \oplus_{j=1}^{m'} Fv_j^*$ and $\mh = \mc \oplus (
\oplus_{i \in I} F\a_i^* )$. Furthermore, we choose and fix a
complement $\mh'' = \oplus_{k=1}^{m'} F\a_{i_k}^*$ of $\mh'$ in
$\mh$, that is, $\mh = \mh' \oplus \mh''$.

Let $\tau \in \Aut(\mb^+)$. Put $\mathfrak{t} = \tau(\mh) \subset
\mb^+$. We would like to show that $\mathfrak{t}$ is a split Cartan
subalgebra (i.e.~a maximal ad-diagonalizable subalgebra) of $\mg$.

We note $\mb^+ = \tau(\mb^+) = \tau(\mh \oplus \mn^+) = \tau(\mh)
\oplus \tau(\mn^+) = \mathfrak{t} \oplus \mn^+$, since
$\tau(\mn^+) = \tau([\mb^+,\mb^+]) = [\tau(\mb^+),\tau(\mb^+)] =
[\mb^+,\mb^+] = \mn^+$. Then, we have $\mn^+ = \oplus_{\gamma \in
\Delta^+(\mathfrak{t})} \mn_\gamma^+(\mathfrak{t})$, where
$\mn_\gamma^+(\mathfrak{t}) = \{ x \in \mn^+\,|\, [t',x] =
\gamma(t')x\ (\forall t' \in \mathfrak{t}) \}$ for each $\gamma
\in \mathfrak{t}^*$, and where $\D^+(\mathfrak{t}) = \{\gamma \in
\mathfrak{t}^*\,|\, \mn_\gamma^+(\mathfrak{t}) \not= 0\}$.
Actually $\D^+(\mathfrak{t})=\D^+\circ\tau^{-1}$, and
$\mn^+_{\gamma}(\mathfrak{t}) = \tau(\mg_\a)$ with
 $\gamma=\a\circ\tau^{-1}$. We define $\gamma_i \in
\mathfrak{t}^*$ by $[t',\tau(e_i)] = \gamma_i(t') \tau(e_i)$ for all
$t' \in \mathfrak{t}$. Then, we see $\gamma_i=\a_i\circ\tau^{-1}$ and
$\D^+(\mathfrak{t}) \subset 
\sum_{i \in I} \Bbb{Z}^+ \gamma_i$.

Put $\mathfrak{p}_i = Ff_i \oplus \mb^+$. Then
$[\mathfrak{t},\mathfrak{p}_i] \subset [\mb^+,\mathfrak{p}_i]
\subset \mathfrak{p}_i$. Hence, in particular, the action of $\mathfrak{t}$
on $\mathfrak{p}_i$ is locally finite. We fix $i$ for a moment.

Let $t' \in \mathfrak{t}$, and write $[t',f_i] = af_i + y$ with $a
\in F$ and $y \in \mb^+$. Since $(Ff_i\,|\,\mn^+) =
(Ff_i\,|\,Fe_i) \not= 0$, we can find $\gamma' \in
\D^+(\mathfrak{t})$ and $z_{\gamma'} \in \mn_{\gamma'}^+(\mathfrak{t})$ such
that $(f_i\,|\,z_{\gamma'}) \not= 0$. Using this $z_{\gamma'}$, we
can obtain
$$a(f_i\,|\,z_{\gamma'}) = (af_i+y\,|\,z_{\gamma'}) = ([t',f_i]\,|\,z_{\gamma'})
= - (f_i\,|\,[t',z_{\gamma'}]) = -\gamma'(t')(f_i\,|\,z_{\gamma'})$$
and $a = -\gamma'(t')$.

Now, we choose $t_0 \in \mathfrak{t}$ such that $-\gamma'(t_0) \not=
\gamma(t_0)$ for all $\gamma \in \Delta_+(\mathfrak{t}) \cup \{ 0 \}$.
Therefore, $-\gamma'(t_0)$ is a single (or multiplicity free) eigenvalue for $t_0$ on $\mathfrak{p}_i$.
Thus, we can find $x_0 \in \mathfrak{p}_i$ such that $[t_0,x_0] = -\gamma'(t_0)x_0$
and $\mathfrak{p}_i = Fx_0 \oplus \mathfrak{t} \oplus \mn^+$
by the local finiteness of $t_0$ on $\mathfrak{p}_i$.
Since $[t_0,[t',x_0]] = [t',[t_0,x_0]] = -\gamma'(t_0)[t',x_0]$
for each $t' \in \mathfrak{t}$,
we see that $[t',x_0]$ must be in $Fx_0$.
Hence, $x_0$ is a common eigenvector for $\mathfrak{t}$.
This means that $\mathfrak{t}$ is diagonalizable on $\mathfrak{p}_i$
under the adjoint action.
In particular, $\mathfrak{p}_i$ is a $\mathfrak{t}^*$-graded
$\mathfrak{t}$-module.

Since $\mg$ is generated by $\mathfrak{p}_0,\mathfrak{p}_1,\ldots,\mathfrak{p}_\ell$,
we obtain that $\mg$ is also a $\mathfrak{t}^*$-graded
$\mathfrak{t}$-module.
Therefore, $\mathfrak{t}$ must be contained in some split Cartan subalgebra
of $\mg$. Comparing both dimensions, we can see that
$\mathfrak{t}$ itself is a split Cartan subalgebra of $\mg$
(cf. \cite{PK}).

Using the conjugacy result for split Cartan subalgebras of $\mg$
(cf.~\cite{PK}), we find an element $g \in G$ such that
$g(\mathfrak{t}) = \mh$. We write $g = u_- w u_+$ with $u_\pm \in
U_\pm$ and $w \in N$. Set $\mathfrak{t}' = u_+(\mathfrak{t})
\subset \mb^+$ and $\mathfrak{t}'' = u_-^{-1}(\mh) \subset \mb^-$.
Then, we see $w(\mathfrak{t}') = \mathfrak{t}''$. This means
$\mathfrak{t}' \subset \mb^+(w):=\mb^+ \cap w^{-1}(\mb^-)$. Since
$\mb^+(w)$ is finite dimensional, we find an element $u$ of the
inner automorphism group, called ${\rm Int}(\mathfrak{b}^+(w))$,
of $\mathfrak{b}^+(w)$, which is corresponding to a certain element $u'
\in U_+$, such that $u(\mathfrak{t}') = u'(\mathfrak{t}') = \mh$
(cf.~\cite{H}, \cite{J}, \cite{MP}; Sect.7, p.641, Ex.7.6). Put
$u'' = u' u_+ \in U_+$. Then, $u''(\mathfrak{t}) = \mh$.
Therefore, we can assume $\tau(\mh) = \mh$ modulo $U_+$.

Since $\tau(\mh) = \mh$, for each $\alpha \in \D_+$ there exists $\b
\in \D_+$ such that $\tau(\mg_\a) = \mg_\b$. Hence, $\tau$ must
permute the $\mg_{\a_i}$. Checking the condition
$[\mg_{\a_i},[\mg_{\a_i},[ \cdots,[\mg_{\a_i},\mg_{\a_j}]\cdots]]] =
0$, the permutation $\sigma$ of $I$, induced by $\tau$, have the
property $a_{\sigma(i),\sigma(j)} = a_{ij}$ for all $i,j$. That is,
$\sigma \in \Aut(A)$. We introduce an automorphism, called $\tilde
\sigma$, of $\mb^+$ defined by $\tilde \sigma(h_0) = h_0,\ \tilde
\sigma(\a_i^*) = \a_{\sigma(i)}^*,\ \tilde \sigma(e_i) =
e_{\sigma(i)}$ for all $h_0 \in \mc$ and $i \in I$. Thus, we can
assume $\tau(\mg_{\a_i}) = \mg_{\a_i}$ for all $i \in I$ modulo
$\Aut(A)$.

Considering the action of $T$ on $\mn^+$,
we can assume $\tau(e_i) = e_i$ for all $i \in I$ modulo $T$.
Then, we obtain
$$\tau(\a_i^*) = \a_i^* + z_i\ (z_i \in \mc),
\quad \tau|_{\mc} \in GL(\mc).$$ Put $\tilde \Gamma = \{ \tilde
\sigma\,|\, \sigma \in \Aut(A) \} \subset \Aut(\mb^+)$. We define
$\Gamma_0(\mb^+) = \{ \tau \in \Aut(\mb^+)\,|\, \tau(\mh) = \mh,\,
\tau(e_i) = e_i\, (i \in I) \}$ and $\Gamma(\mb^+) = \langle \Gamma_0(\mb^+),\,
\tilde \Gamma \rangle$. Then $\Gamma(\mb^+)/\Gamma_0(\mb^+) \simeq \Aut(A)$, and
we see $\Gamma_0(\mb^+) = GL(\mc) \ltimes \mc^{l+1}$ and
$\Gamma(\mb^+) = (\tilde \Gamma \times GL(\mc))
\ltimes \mc^{l+1}$. Therefore, we just established the
following.

\begin{thm} Let $A$ be an indecomposable symmetrizable GCM. Then:
$$\Aut(\mb^+) = \Gamma(\mb^+) \ltimes B_+ =
\Gamma(\mb^+) \ltimes (T \ltimes U_+)\ .$$
In particular, if $A$ is nonsingular, then we have
$$\Aut(\mb^+) = \Aut(A) \ltimes B_+ = \Aut(A) \ltimes (T \ltimes U_+)\ .$$
\end{thm}

Here,
we make a remark on the automorphism group $\Aut(\mg)$ of $\mg$
to compare to $\Aut(\mb^+)$,
which is already discussed in \cite{Kac},\cite{MP},\cite{PK}.
For each $\sigma \in \Aut(A)$, we define
$\sigma^\star \in \Aut(\mg)$ by
$\sigma^\star(\a_i^\vee) = \a_{\sigma(i)}^\vee,\
\sigma^\star(\a_{i_k}^*) = \a_{\sigma(i_k)}^*,\
\sigma^\star(e_i) = e_{\sigma(i)},\ \sigma^\star(f_i) = f_{\sigma(i)}$
for all $i \in I$ and $1 \leq k \leq m'$ (cf.~\cite{Kac},\cite{MP}).
Let $\Gamma(\mg) = \{ \rho \in \Aut(\mg)\, |\,
\exists \sigma \in \Aut(A)\ {\rm s.t.}\
\rho(\mh) = \mh,\, \rho(e_i) = e_{\sigma(i)},\,
\rho(f_i) = f_{\sigma(i)}\ (i \in I) \}$ and
$\Gamma_0(\mg) = \{ \rho \in \Gamma(\mg)\, |\,
\rho(e_i) = e_i,\ \rho(f_i) = f_i\ (i \in I) \}$.
Then,
$\mc^{m'} \simeq \Gamma_0(\mg) \subset \Gamma(\mg)$ and
$\Gamma(\mg)/\Gamma_0(\mg) \simeq \Aut(A)$ given by
$(\sigma^\star\, {\rm mod}\, \Gamma_0(\mg)) \leftrightarrow \sigma$.
By the conjugacy theorem for split Cartan subalgebras of $\mg$
(cf.~\cite{Kac},\cite{MP},\cite{PK}), we obtain the following result.

\begin{thm} {\rm (cf.~\cite{PK})}
Let $A$ be an indecomposable  symmetrizable GCM. Then:
$$\Aut(\mg) = \langle \Gamma(\mg), \omega, G \rangle\ ,$$
especially
$\Aut(\mg) = (\langle \omega \rangle \ltimes \Gamma(\mg)) \ltimes G$
if $A$ is of infinite type.
Furthermore, if $A$ is nonsingular, then we have
$$\Aut(\mg) = \langle \Aut(A), \omega, G \rangle\ ,$$
more precisely
$\Aut(\mg) = \Aut(A) \ltimes G$ if $A$ is of finite type,
and
$\Aut(\mg) = (\Aut(A) \times \langle \omega \rangle ) \ltimes G$
if $A$ is of infinite type.
\end{thm}

\section{\textbf{ Derivations of $\mn^+$ and $\tn^+$ }}

In this section we need not assume that $A$ is symmetrizable. We
will determine $\Der(\mn^+)$ and $\Der(\tn^+)$ in general, which
is useful for us to study the automorphisms of the Lie algebra
$\mn^+$. About 30 years ago, Moody conjectured in \cite{M}:
\textquoteleft \textquoteleft It is conceivable that $\Der(\tn^+)
= \ad(\tb^+)|_{\tn^+}$ if $\dim \tn^+$ is infinite.'' When $\dim
\tn^+$ is finite, $\Der(\tn^+)=\Der(\mn^+)$ was already known
(cf.~\cite{Ko},\cite{LL2}). In \cite{F2} 1993,
$\Der(\tn^+)=\Der(\mn^+)$ was established for affine type, which
particularly shows $\Der(\tn^+) \ne \ad(\tb^+)|_{\tn^+}$ (also see
\cite{F1}). Here, we will use elementary methods to determine
$\Der(\mn^+)$ and $\Der(\tn^+)$. Consequently, we recover
Kostant's and Fialowski's results (cf.~\cite{Ko}, \cite{F2}) and
completely solve Moody's conjecture.

We shall mainly discuss the algebra $\mn^+$. The same arguments
are valid for $\tn^+$ also, that is, we have to keep in mind that
$\mn^+$ represents $\mn^+$ and $\tn^+$.

We denote by $End(\mn^+)$ the space of all $F$-linear maps from
$\mn^+$ to itself. For each $\b \in Q$, we set $End(\mn^+)_\b =
\left\{ \phi' \in End(\mn^+)\, |\, \phi'(\mg_\a) \subset \mg_{\a +
\b}\ (\forall\a \in \D_+) \right\}$. For $\phi = \prod \phi_\b \in
\displaystyle{\prod_{\b \in Q} End(\mn^+)_\b}$ and $\a \in \D^+$,
we define $\Omega(\phi,\alpha) = \{ \b \in Q \mid \phi_\b(\mg_\a)
\ne 0 \}$. Then, we have
$$End(\mn^+) = \left\{ \left. \phi \in
\prod_{\b \in Q} End(\mn^+)_\b \right| \sharp\, \Omega(\phi,\alpha) < \infty\
(\forall \a \in \D_+) \right\} \subset \prod_{\b \in Q} End(\mn^+)_\b\ .$$
For $d = \prod d_\b \in End(\mn^+) \subset \displaystyle{\prod_{\b \in Q} End(\mn^+)_\b}$,
we easily see
$d \in \Der(\mn^+) \Leftrightarrow d_\b \in \Der(\mn^+)$
for all $\b \in Q$.
Put $\Omega(d) = \cup_{i \in  I}\ \Omega(d,\a_i)$ for each $d \in
End(\mn^+)$. If $d \in \Der(\mn^+)$, then $d = \sum_{\b \in
\Omega(d)}\, d_\b$, since $\mn^+$ is generated by $e_1,\cdots,e_l$
and $\Omega(d)$ is finite. Therefore, we see
$$\Der(\mn^+) = \oplus_{\b \in Q}\ \Der(\mn^+)_\b\ ,$$
where $\Der(\mn^+)_\b = \Der(\mn^+) \cap End(\mn^+)_\b$.

For our purpose to determine $\Der(\mn^+)$, we may assume $d = d_\b
\in \Der(\mn^+)_\b$, i.e.,   there exists $\b\in Q$ such that
$$d(e_i)=e_{\b+\a_i},\quad \forall\,\,i\in I,\eqno(3.1)
$$
where  $e_{\b+\a_i}\in \mn^+_{\b+\a_i}$ with $e_{\b+\a_i} = 0$ if
$\b+\a_i \notin \D_+$. Note that, in general, $e_i\ne e_{\a_i}$ in
our notation. Let $I_1=\{i\in I\,|\, d(e_i)\ne0\}$ and $I_0=\{i\in
I\,|\, d(e_i)=0\}$.

If $\b=0$ in (3.1) and $d(e_i) = e_{\a_i}=c_ie_i$ for $i \in I$,
then $d = \ad h$ with $h = \sum_{i \in I}\, c_i\a_i^*$.

Next we assume that
$d\ne0$, $\b\ne0$ and $A$ is not of type $A_1$ (because it is
trivial).

\begin{lem} If
$[f_i,e_{\b+\a_i}]=0$ for all $i\in I$, then $A$ is of finite type, $I_1=\{i_1\}$ and
$\b=s_{i_1}(\theta)-\a_{i_1}$, where $s_{i_1}\in W$ is the
reflection determined by $\a_{i_1}$ and $\theta$ is the highest
root of $\mg = \mg(A)$. \label{Le:1}
\end{lem}

\noindent {\it Proof.} We will prove this lemma in several steps.

\noindent {\bf Claim 1}. {\it $\b(\a_{i_1}^\vee)\le -2$ for all
$i_1\in I_1$.}

\noindent Since $[f_i,e_{\b+\a_i}]=0$ for all $i\in I$, then
$(\b+\a_{i_1})(\a_{i_1}^\vee)\le 0$, i.e., $\b(\a_{i_1}^\vee)\le
-2$ for all $i_1\in I_1$.

\noindent {\bf Claim 2}. $|I_1|=1$.

\noindent Suppose $i_1, i_2\in I_1$ are different, and let
$r=-\a_{i_2}(\a_{i_1}^\vee)$, $r_1=-(\b+\a_{i_2})(\a_{i_1}^\vee)$.
Then $r_1\ge r+2$, $(\ad e_{i_1})^{r+1}e_{\b+\a_{i_2}}\ne0$ and
further $(\ad f_{i_1})^{r+1}(\ad
e_{i_1})^{r+1}e_{\b+\a_{i_2}}\ne0$ (using $sl_2$-module structure).
We deduce that
$$0=d((\ad
e_{i_1})^{1+r}e_{i_2})=(\ad
e_{i_1})^{1+r}e_{\b+\a_{i_2}}+\sum_{s=0}^r(\ad e_{i_1})^s(\ad
e_{\b+\a_{i_1}})(\ad e_{i_1})^{r-s} e_{i_2}.$$ Applying $(\ad
f_{i_1})^{r+1}$ to the above equation we obtain  $(\ad
f_{i_1})^{r+1}(\ad e_{i_1})^{r+1}e_{\b+\a_{i_2}}=0$,
a contradiction.
(Use $[f_{i_1},e_{\b+\a_{i_1}}]=[f_{i_1},e_{i_2}]=0$
to show
$(\ad f_{i_1})^{r+1}(\ad e_{i_1})^s(\ad e_{\b+\a_{i_1}})(\ad e_{i_1})^{r-s} e_{i_2} = 0.$)
This proves Claim 2.

Suppose $I_1=\{i_1\}$. Let $\b_1=\b+\a_{i_1}$,
$r_1=-\b_1(\a_{i_1}^\vee)$, $Y=(\ad e_{i_1})^{r_1}e_{\b_1}$. Then
$r_1\ge0$, $Y\ne0$ and $[ e_{i_1}, Y]=0$. Fix $i\in I_0$, let
$r=-\a_i(\a_{i_1}^\vee)$.

\noindent {\bf Claim 3}. $[ e_{i}, Y]=0$.

\noindent We shall break the proof into two cases.

{\bf Case 1}: $r_1=0$.

If $r=0$, we know that $0=d([e_i, e_{i_1}])=[e_i, e_{\b_1}]=[
e_{i}, Y]$. If $r>0$, we deduce that
$$0=d((\ad
e_{i_1})^{1+r}e_{i})=\sum_{s=0}^r(\ad e_{i_1})^s(\ad e_{\b_1})(\ad
e_{i_1})^{r-s} e_{i} =(r+1) (\ad e_{i_1})^r[e_{\b_1},  e_{i}].$$
If $[e_{\b_1},  e_{i}]\ne0$, then $(\ad e_{i_1})^r[e_{\b_1},
e_{i}]=0$ which is impossible since
$-(\b_1+\a_i)(\a_{i_1}^\vee)=r$. Thus $[ e_{i}, e_{\b_1}]=0$,
i.e., $[ e_{i}, Y]=0$.

{\bf Case 2}: $r_1>0$.

If $r=0$, we know $[ e_{i}, e_{i_1}]=0$. Then $0=
d([e_{i},e_{i_1}])=[e_{i},d(e_{i_1})]=[ e_{i},e_{\b_1}]$, and
further $[ e_{i}, Y]=[e_i,(\ad e_{i_1})^{r_1}e_{\b_1}]= (\ad
e_{i_1})^{r_1}[e_i,e_{\b_1}]=0$.  Now suppose
$r=-\a_i(\a_{i_1}^\vee)>0$.

If $r_1\ge r$, we compute
$$0=d((\ad
e_{i_1})^{1+r}e_{i})=\sum_{s=0}^r(\ad e_{i_1})^s(\ad e_{\b_1})(\ad
e_{i_1})^{r-s} e_{i}=X.\eqno(3.2)$$

We use the notation $[x_1,x_2,x_3,\cdots,x_{k-1},x_k]$ for
$[x_1,[x_2,[x_3,[\cdots,[x_{k-1},x_k]\cdots]]]]$. For $r' \geq s
\geq 0$, we put $P(r',s ) = (\ad e_{i_1})^{s }(\ad e_{\b_1})(\ad
e_{i_1})^{r'-s }e_i$. Note that $[f_{i_1},e_{\b_1}] = 0$. Then, for
$r \geq r' \ge 0$, we have
$$\begin{array}{ll}
& [f_{i_1},P(r',s )]\
=\ [f_{i_1},(\ad e_{i_1})^{s }(\ad e_{\b_1})(\ad e_{i_1})^{r'-s }e_i]\\
&\\
= & (\ad f_{i_1}) [ \underbrace{e_{i_1},\cdots,e_{i_1}}_{s }
,e_{\b_1}, \underbrace{e_{i_1},\cdots,e_{i_1}}_{r'-s }
,e_i]\\
= & \sum [ \underbrace{e_{i_1},\cdots,e_{i_1}}_{s } ,e_{\b_1},
\underbrace{e_{i_1},\cdots,-\a_{i_1}^\vee,\cdots,e_{i_1}}_{r'-s }
,e_i]\\
&\quad + \sum [
\underbrace{e_{i_1},\cdots,-\a_{i_1}^\vee,\cdots,e_{i_1}}_{s }
,e_{\b_1}, \underbrace{e_{i_1},\cdots,e_{i_1}}_{r'-s } ,e_i]\\
 = &
\sum\limits_{t=0}^{r'-s -1} (\a_i+t\a_{i_1})(-\a_{i_1}^\vee)
(\ad e_{i_1})^{s }(\ad e_{\b_1})(\ad e_{i_1})^{r'-s -1}e_i\\
&\\ [-0.2cm] & \quad + \sum\limits_{u=0}^{s -1} (\a_i+(r'-s
)\a_{i_1}+\b_1+u\a_{i_1})(-\a_{i_1}^\vee)
(\ad e_{i_1})^{s -1}(\ad e_{\b_1})(\ad e_{i_1})^{r'-s }e_i
\end{array}$$
$$
\begin{array}{ll}
&\\
 = & (\sum\limits_{t=0}^{r'-s -1} (r-2t))\cdot
(\ad e_{i_1})^{s }(\ad e_{\b_1})(\ad e_{i_1})^{r'-s -1}e_i\\
&\\ [-0.2cm] & \quad +\, (\sum\limits_{u=0}^{s -1} (r-2(r'-s
)+r_1-2u))\cdot
(\ad e_{i_1})^{s -1}(\ad e_{\b_1})(\ad e_{i_1})^{r'-s }e_i\\
&\\ [-0.2cm] = & (r-r'+s +1)(r'-s )\cdot
({\rm ad}\ e_{i_1})^{s }({\rm ad}\ e_{\b_1})({\rm ad}\ e_{i_1})^{r'-s -1}e_i\\
&\\ [-0.2cm] & \quad +\, (r-2r'+s +r_1+1)s \cdot
(\ad e_{i_1})^{s -1}(\ad e_{\b_1})(\ad e_{i_1})^{r'-s }e_i\\
&\\ [-0.2cm] = & (r-r'+s +1)(r'-s )\cdot P(r'-1,s ) + (r-2r'+s
+r_1+1)s \cdot P(r'-1,s -1)\ ,
\end{array}
$$
where $(r-r'+s +1)(r'-s )$ is  positive at most one is zero, and
$(r-2r'+s +r_1+1)s $ is positive at most one is zero since $r_1 \geq
r$. Therefore,
$(\ad f_{i_1})^rX=(\ad f_{i_1})^r \sum_{s=0}^r\, P(r,s)
= aP(0,0) = 0$
with some positive integer $a$, yielding
$$P(0,0)=[e_{\b_1},e_i]=0.\eqno(3.3)$$
In particular, $P(1,1) = 0$.
Using the fact
$(\ad f_{i_1})^{r-1}X
= bP(1,1) + cP(1,0) = cP(1,0)$
for some positive integers $b$ and $c$,
we obtain that
$P(1,0)=0$. In particular, $P(2,2) = P(2,1) = 0$.
Using the fact
$(\ad f_{i_1})^{r-2}X
= a'P(2,2) + b'P(2,1) + c'P(2,0) = c'P(2,0)$
for some positive integers $a',b',c'$,
we obtain that
$P(2,0)=0$. In particular, $P(3,3) = P(3,2) = P(3,1) = 0$.
Continuing this, we have
$P(k,j) = 0\quad \mbox{for all}\ 0 \leq j \leq k \leq r$.
Now let $Q(k,j) = [(\ad e_{i_1})^j e_{\b_1},(\ad e_{i_1})^{k-j} e_i]$ for $0 \leq j \leq k \leq r$.
By the definition, we note $Q(k,0) = P(k,0) = 0$.
Since $Q(k,j) = [e_{i_1},Q(k-1,j-1)] - Q(k,j-1)$ with $j > 0$,
we can inductively obtain $Q(k,j) = 0$ for all $0 \leq j \leq k \leq r$.
In particular, we obtain that
$$[e_i,(\ad e_{i_1})^ke_{\b_1}]=-Q(k,k)=0, \,\,\forall\,\, k\le r.\eqno(3.4)$$
If $r_1>r$, using the above established formulas, we consider
$Q(r+1,0) = [e_{\b_1},(\ad e_{i_1})^{r+1}e_{i}]=0$ to obtain that
(3.4) with $k=r+1$. In this manner, repeating this process, we can
obtain that $[e_i,Y]=0$.

Now suppose $r>r_1>0$. We first recall a couple of formulas on
$sl_2$-module structure: \begin{align*} & (\ad f_{i_1})(\ad
e_{i_1})^ke_{\b_1}=k(r_1+1-k)(\ad
e_{i_1})^{k-1}e_{\b_1},\,\forall \, 0\le k\le r_1,\tag{3.5}\\
& (\ad f_{i_1})(\ad e_{i_1})^ke_{i}=k(r+1-k)(\ad
e_{i_1})^{k-1}e_{i},\,\forall \, 0\le k\le r.\tag{3.6}
\end{align*}

We know that
\begin{align*}
0 = d((\ad e_{i_1})^{1+r}e_{i}) & = \sum_{s=0}^r(\ad
e_{i_1})^s(\ad e_{\b_1})(\ad e_{i_1})^{r-s} e_{i}\end{align*}
\begin{align*}
& = \sum_{s=0}^r(\ad e_{i_1})^s([e_{\b_1},(\ad e_{i_1})^{r-s} e_{i}]\\
& = \sum_{s=0}^r\sum_{j=0}^s {s\choose{j}} [(\ad e_{i_1})^j e_{\b_1},
(\ad e_{i_1})^{s-j} (\ad
e_{i_1})^{r-s} e_{i}]\\
& = \sum_{s=0}^r\sum_{j=0}^s {s\choose{j}} [(\ad e_{i_1})^j
e_{\b_1}, (\ad e_{i_1})^{r-j} e_{i}]\\ & =
\sum_{s=0}^r\sum_{j=0}^{\min\{s, r_1\}} {s\choose{j}} [(\ad
e_{i_1})^j e_{\b_1},
(\ad e_{i_1})^{r-j} e_{i}]\\
& = \sum_{j=0}^{r_1}\left( \sum_{s=j}^r {s\choose{j}}\right) [(\ad e_{i_1})^j e_{\b_1},
(\ad e_{i_1})^{r-j} e_{i}]\\
& = \sum_{j=0}^{r_1} {r+1\choose{j+1}} [(\ad e_{i_1})^j e_{\b_1},
(\ad e_{i_1})^{r-j} e_{i}] =  X.\tag{3.7}
\end{align*}
Here, we notice $\{\, (s,j) \mid  0 \leq s \leq r,\ 0 \leq j \leq
{\rm min}\{s,r_1\}\, \} = \{\, (s,j) \mid 0 \leq j \leq r_1,\ j
\leq s \leq r\, \}$. Using (3.5) and (3.6) we know that $0=(\ad
f_{i_1})^rX=a[e_{\b_1},e_i]$ where $a$ is a nonzero constant. Thus
$$[e_{\b_1},e_i]=0.\eqno(3.8)$$
From $0=(\ad e_{i_1})^{r+r_1}[e_{\b_1},e_i]$, we deduce that
$$[(\ad e_{i_1})^{r_1}e_{\b_1},(\ad e_{i_1})^{r}e_i]=0.\eqno(3.9)$$
By induction on $k: 0\le k\le r_1$, we will prove
\begin{align*}
&[(\ad e_{i_1})^{j}e_{\b_1},(\ad e_{i_1})^{k-j}e_i]=0,\quad 0\le j\le k
,\quad {\rm and} \tag{3.10}\\
&[(\ad e_{i_1})^{r_1-k+j}e_{\b_1},(\ad e_{i_1})^{r-j}e_i]=0,\quad
0\le j\le k.\tag{3.11}
\end{align*}
The formulas (3.8) and (3.9) ensure (3.10) and (3.11)  for $ k=0$.
Suppose that (3.10) and (3.11) hold for $0\le k\le k_0$ where
$k_0\le r_1-1$. Let us consider (3.10) and (3.11) for $k=k_0+1$.
By our induction, we have
$$\begin{array}{l}
\left[ (\ad e_{i_1})^{j}e_{\b_1},(\ad e_{i_1})^{r_1+r-k_0-j}e_i
\right] = 0
\\
\\
 \hskip 4cm\mbox{by}\quad \left\{
\begin{array}{l}
(\ad e_{i_1})^{r_1+r-k_0-j}e_i = 0\quad \mbox{if}\quad r_1-k_0-j > 0,\\
(3.11)\quad 0 \leq j'=j-r_1+k_0 \leq k = k_0\quad \mbox{otherwise};
\end{array}
\right.
\\
\\
\left[ (\ad e_{i_1})^{j}e_{\b_1},(\ad e_{i_1})^{k_0-j}e_i
\right] = 0 \quad \mbox{by}\quad (3.10)\quad 0 \leq j \leq k =
k_0\ .
\end{array}$$
Hence, using
$$\left[ f_{i_1},\left[ (\ad e_{i_1})^{j}e_{\b_1},
(\ad e_{i_1})^{r_1+r-k_0-j}e_i \right] \right] = 0
\quad \mbox{and}\quad
\left[ e_{i_1},\left[ (\ad e_{i_1})^{j}e_{\b_1},
(\ad e_{i_1})^{k_0-j}e_i \right] \right] = 0\ ,$$
we deduce that
\begin{align*}
& [(\ad e_{i_1})^{j-1}e_{\b_1},(\ad e_{i_1})^{r_1+r-k_0-j}e_i]\\
& \qquad =\frac{-(r_1+r-k_0-j)(k_0+j+1-r_1)}{j(r_1+1-j)}[(\ad
e_{i_1})^{j}e_{\b_1},(\ad e_{j})^{r_1+r-k_0-j-1}e_i]\ ,\tag{3.12}\\
& [(\ad e_{i_1})^{j+1}e_{\b_1},(\ad e_{i_1})^{k_0-j}e_i]=-[(\ad
e_{i_1})^{j}e_{\b_1},(\ad e_{i_1})^{k_0-j+1}e_i]\ .\tag{3.13}
\end{align*}
From (3.12) we recursively obtain that
\begin{align*}
& [(\ad
e_{i_1})^{j}e_{\b_1},(\ad e_{i_1})^{r_1+r-k_0-j-1}e_i]\\
& \qquad
=(-1)^{r_1-j}\frac{{r+r_1-k_0-j-1\choose{r_1-j}}{k_0+1\choose{r_1-j}}}{{r_1\choose{j}}}[(\ad
e_{i_1})^{r_1}e_{\b_1},(\ad e_{j})^{r-k_0-1}e_i]\ .\tag{3.14}
\end{align*}
Using (3.7) and (3.14) we deduce that
\begin{align*}
0 & = (\ad e_{i_1})^{r_1-k_0-1}X\\
& = \sum_{j=0}^{r_1} {r+1\choose{j+1}} \sum_{s=0}^{r_1-k_0-1} {r_1-k_0-1\choose{s}}[(\ad
e_{i_1})^{s+j} e_{\b_1}, (\ad e_{i_1})^{r+r_1-k_0-s-j-1} e_{i}]\\
& = \sum_{k=r_1-k_0-1}^{r_1} \,\, \sum_{s+j=k}{r+1\choose{j+1}} {r_1-k_0-1\choose{s}}[(\ad
e_{i_1})^{k} e_{\b_1}, (\ad e_{i_1})^{r+r_1-k_0-k-1} e_{i}]\\
& = \sum_{k=r_1-k_0-1}^{r_1} {r+r_1-k_0\choose{k+1}}  [(\ad
e_{i_1})^k e_{\b_1}, (\ad e_{i_1})^{r+r_1-k_0-k-1} e_{i}]\\
& = \sum_{k=r_1-k_0-1}^{r_1}
\frac{(-1)^{r_1-k}{r+r_1-k_0\choose{k+1}}{r+r_1-k_0-k-1\choose{r_1-k}}
{k_0+1\choose{r_1-k}}}{{r_1\choose{k}}} [(\ad e_{i_1})^{r_1}
e_{\b_1}, (\ad e_{i_1})^{r-k_0-1} e_{i}]\ ,\tag{3.15}
\end{align*}
where we have used the combinatorics formula $
\sum\limits_{s+j=k}{r+1\choose{j+1}} {r_1-k_0-1\choose{s}}=
{r+r_1-k_0\choose{k+1}} $. We compute the coefficient on the
right-hand side of (3.15):
\begin{align*}
& \sum_{k=r_1-k_0-1}^{r_1}
\frac{(-1)^{r_1-k}{r+r_1-k_0\choose{k+1}}{r+r_1-k_0-k-1\choose{r_1-k}}
{k_0+1\choose{r_1-k}}}{{r_1\choose{k}}}\end{align*}
\begin{align*} &= \sum_{k=0}^{k_0+1}
\frac{(-1)^{r_1-k}{k_0+1\choose{k}}{r+r_1-k_0\choose{r-k}}{r-k\choose{k_0+1-k}}}
{{r_1\choose{k_0+1-k}}}\\
&= \frac{(-1)^{r_1}(r_1+1)(r_1+2)...(r_1+r-k_0)}{(r-k_0-1)!}
\sum_{k=0}^{k_0+1}(-1)^k{k_0+1\choose{k}}\frac 1{r_1-k_0+k}\\
&= \frac{(-1)^{r_1}(r_1+1)(r_1+2)...(r_1+r-k_0)}{(r-k_0-1)!}
\sum_{k=0}^{k_0+1}(-1)^k{k_0+1\choose{k}}\frac{x^{r_1-k_0+k}}
{r_1-k_0+k}\Big|_{x=1}\\
&= \frac{(-1)^{r_1}(r_1+1)(r_1+2)...(r_1+r-k_0)}{(r-k_0-1)!}\int_0^1
\sum_{k=0}^{k_0+1}(-1)^k{k_0+1\choose{k}}{x^{r_1-k_0+k-1}}dx\\
&=
\frac{(-1)^{r_1}(r_1+1)(r_1+2)...(r_1+r-k_0)}{(r-k_0-1)!}\int_0^1
(1-x)^{k_0+1} x^{r_1-k_0-1}dx\ne0.\tag{3.16}
\end{align*}
Applying this to (3.15) we see that $[(\ad e_{i_1})^{r_1}
e_{\b_1}, (\ad e_{i_1})^{r-k_0-1} e_{i}]=0$. Using (3.12) we
obtain  (3.11) for $k=k_0+1$.

Now let us deduce (3.10)  for $k=k_0+1$.

Suppose  $[(\ad e_{i_1})^{k_0+1}e_{\b_1},(\ad e_{i_1})e_i]\ne0$.
From $(\b_1+\a_i+(k_0+1)\a_{i_1})(\a_{i_1}^\vee)=-r-r_1+2k_0+2<0$,
we know that $$(\ad e_{i_1})^{r+r_1-2k_0-2}[(\ad
e_{i_1})^{k_0+1}e_{\b_1},(\ad e_{i_1}) e_i]\ne0.$$ Using (3.5) and
(3.6) we know that
$$\begin{array}{lll}
0& \ne & (\ad e_{i_1})^{r+r_1-2k_0-2}[(\ad
e_{i_1})^{k_0+1}e_{\b_1},(\ad e_{i_1}) e_i]\\
&&\\
& \in & \sum_{j=0}^{r_1}\, F[(\ad e_{i_1})^{r_1-j}e_{\b_1},(\ad
e_{i_1})^{r-k_0-1+j}e_i]\ =\ \{ 0 \}\ ,
\end{array}$$
which is a contradiction (used   (3.11) for $k=k_0+1$ at the last
step). Therefore, we see $[(\ad e_{i_1})^{k_0+1}e_{\b_1},(\ad
e_{i_1})e_i]=0$. From (3.13) we get (3.10) for $k=k_0+1$.

Therefore, (3.10) and (3.11) hold for $k\le r_1$. Taking $k=j=r_1$
in (3.10) we see that $[ e_{i}, Y]=0$. Thus $[ e_{i}, Y]=0$ for
all $i\in I$, consequently $A$ is of finite type and
$\b=s_{i_1}(\theta)-\a_{i_1}$. \qed

Lemma 3.1 solved our problem in the case when $\b\notin\Delta_+$.
Next we consider the case of $\b\in\Delta_+$.

\begin{lem} Suppose  $\b\in\Delta_+$. Then
$e_{\b+\a_{i_1}}\in[e_{i_1}, \mn^+_\b]$ for all $i_1\in I_1$.
\end{lem}

\noindent {\it Proof.} To the contrary, suppose there is $i_1\in
I_1$ such that $e_{\b+\a_{i_1}}\notin[e_{i_1}, \mn^+_\b]$. Then
$A$ is not of finite type, and $\b_1=\b+\a_{i_1}\in\Delta^{im}$.

By subtracting an element in $\ad(\mn^+_{\b})$ from $d$, we may
assume that $[f_{i_1}, e_{\b+\a_{i_1}}]=0$. Then
$(\b+\a_{i_1})(\a_{i_1}^\vee)\le0$. Using the same argument as in
the proof of Claim 2 in Lemma 3.1, we see that $|I_1|=1$.

Now all assumptions in Lemma 3.1 are satisfied by our new $d$. By
the same Lemma, we deduce that $A$ is of finite type which is
impossible. Consequently, Lemma 3.2 follows. \qed

From now on, we will use the construction of Kac-Moody algebras of
affine type $A=X_{N}^{(r)}$ $(r=1,2,3$) introduced in Sections 7.4
and 8.3 at \cite{Kac}. We will try to adopt notations from
\cite{Kac}, but we have to make some slight modifications. In this
paper, $\dot \mg$ denotes a finite dimensional simple Lie algebra
of type $X_N$ with a split Cartan subalgebra $\dot \mh$, and
$\mg=\mg(A)$ will be our affine Lie algebra obtained form $\dot
\mg$. (Our $\dot \mg$ is corresponding to $\mg$ in Section 8 of
\cite{Kac}, and we use $\dot \mg$ instead of $\mg$ here, since
$\mg$ in this paper has a different meaning in many places.) Let
$\dot \D$ be the root system of $\dot \mg$ with respect to $\dot
\mh$, and fix a base $\dot \Pi$ of $\dot \D$, and Chevalley
generators $ E'_i, F'_i, H'_i=[E'_i, F'_i] $ $(i=1,2,\cdots N)$ of
$\dot\mg$. We have the decomposition $\dot \mg=\oplus_{i\in\Z_r}\,
\dot \mg_{\bar i}$ , which is graded by a fixed Dynkin diagram
automorphism $\bar \mu$ of order $r$, satisfying $\dot
\mh=\oplus_{i\in\Z_r}\, \dot \mh_{\bar i}$ with $\dot \mh_{\bar i}
= \dot \mh \cap \dot \mg_{\bar i}$, and we choose Chevalley
generators $E_{i}, F_{i},H_i=[E_i,F_i]$ for $\dot \mg_{\bar 0}$
with $0 \leq i \leq l-1$ if $A=A_{2l}^{(2)}$, and $1 \leq i \leq
l$ otherwise. Then, we can also have the Chevalley generators
$e_i, f_i, \a^\vee_i=[e_i, f_i]$ ($i\in I$) for the affine
Kac-Moody algebra $\mg(A)$. Note  that $\epsilon=l$ for
$A=A_{2l}^{(2)}$ and $\epsilon=0$ for all other affine types. Let
$\delta \in \D^+$ be the smallest positive imaginary root. With
respect to the degree of the variable $t$ in the construction of
affine Lie algebras, we have $\deg(e_\epsilon)=1$ and
$\deg(e_i)=0$ if $i\ne \epsilon$.

\noindent {\bf Remark 1}.\ In the construction of the affine Lie
algebra of type $D_4^{(3)}$, it seems we require that $x^3-1$
splits in $F$. Actually this is not necessary. If $F$ does not
contain a cubic root $\eta$ of $1$, then first we put $F' =
F(\eta)$ and we construct everything over $F'$, that is, $\dot
\mg_{F'}$ and $\mg(A)_{F'}$ etc. (We may naturally take
$\eta=e^{2\pi \sqrt{-1}/3}$ if $F \subset \C$.) Next, we choose
standard Chevalley generators of $\mg(A)_{F'}$, which generate a
subalgebra over $F$, called $\mathfrak{L}_F$, of the loop algebra
corresponding to $[\mg(A)_{F'},\mg(A)_{F'}]/\mc_{F'}$. Then,
$\mathfrak{L}_F$ is a loop algebra associated to our original
affine Lie algebra $\mg(A) = \mg(A)_F$, that is, $\mathfrak{L}_F$
is isomorphic to $[\mg(A)_F,\mg(A)_F]/\mc$, and $\mg(A)_F$ is
isomorphic to (or identified with)
$$\mathfrak{L}_F \oplus F{\boldsymbol z} \oplus Ft\frac{d}{dt}$$
where $\mc = F{\boldsymbol z}$ is the center of $\mg(A)_F$ and contained in $[\mg(A)_F,\mg(A)_F]$.

\begin{lem} If  $d\ (\ne 0) \in \Der(\mn^+)_\b$ with $\b\in\Delta_+$, then

(a). {\rm $d\in \ad(\mn^+)$}, or

(b). $A$ is of affine type $X_{N}^{(r)}$, $\b=kr\delta$ $(k\in\N)$
and there exists $e_\b\in\mn^+_\b$ such that $(d-{\rm ad}
e_\b)(e_i)=\delta_{i,\epsilon}e_{\b+\a_\epsilon}$ for all $i\in I$
where $e_{\b+\a_\epsilon}\in\mn^+_{{\b+\a_\epsilon}}$.

\end{lem}

\noindent {\it Proof.} From Lemma 3.2 we know that
$e_{\b+\a_{i}}\in[e_{i}, \mn^+_{\b}]$ for all $i\in I$. We may
assume that $d(e_i)=[e_{\b}^{(i)}, e_i]$ for all $i\in I$ where
$e_{\b}^{(i)}\in \mg_{\b}$.

If there exists $i_1\in I$ such that $\b(\a_{i_1}^\vee)<0$, by
subtracting $\ad (e_{\b}^{(i_1)})$ from $d$ we may assume that
$e_{\b}^{(i_1)}=0$, i.e., $d(e_{i_1})=0$.   For any $i\ne i_1$, let
$r=-\a_i(\a_{i_1}^\vee)$. Then
$$0=d((\ad
e_{i_1})^{1+r}e_{i})=(\ad e_{i_1})^{1+r}[ e_\b^{(i)},e_{i}].$$ If
$[ e_\b^{(i)},e_{i}]\ne0$, since $(\b+\a_{i})(\a_{i_1}^\vee)<-r$
we know that $(\ad e_{i_1})^{1+r}[ e_\b^{(i)},e_{i}]\ne0$, a
contradiction. Thus $[ e_\b^{(i)},e_{i}]=0$, i.e., $d=0$.

Therefore, as the remaining case we should consider, we have
$\b(\a_i^\vee)\ge 0$ for all $i\in I$. If $\b(\a_i^\vee)> 0$
for at least one $i$, then $A$ is of finite type, and $\b=\theta$
(the highest root) in
which case $d=0$, or $\b=\theta_1$ the highest short root if $A$ is
one of $B_{l}$, $C_{l}$, $F_4$, and $G_2$. By examining the highest
short root for each of $B_{l}$, $C_{l}$, $F_4$, and $G_2$, we see
that there exists exactly one $i_0\in I$ such that
$[e_{i_0},\mg_{\theta_1}]\ne0$ and $[e_{i},\mg_{\theta_1}]=0$ for
all other $i\ne i_0$. (For more details, please refer to Section 12
in \cite{H}. See also \cite{Bo}.) More precisely, $i_0=l$ for $B_{l}$, $i_0=1$ for
$C_{l}$, $i_0=3$ for $F_4$, and $i_0=1$ for $G_2$. Thus
$d=\ad(e_{\theta_1})$ for some $e_{\theta_1}\in \mg_{\theta_1}$.

Now we consider the only case $\b(\a_{i_1}^\vee)= 0$ for all $i\in
I$. Then $A$ is of affine type $X_{N}^{(r)}$ as in TABLE Aff 1,
TABLE Aff 2 and TABLE Aff 3 at \cite{Kac}. If $\b=kr\d$, we can
easily find $e_\b\in\mg(A)_\b=t^{kr}\otimes \dot \mh_{\bar0}$ such that
$\ad e_\b(e_i)=d(e_i)$ for all $i\ne\e$. It is not hard to see that
$(d-\ad e_\b)(e_i)=\delta_{i,\epsilon}e_{\b+\a_\epsilon}$ for all
$i\in I$ determines a derivation of $\mn^+$.

If $A=A_{2l}^{(2)}$ and $\b=(2k+1)\d$, because $\dim\mg(A)_\b=l$
we can  find $e_\b\in\mg(A)_\b$ such that $[e_\b,e_i]=d(e_i)$ for
all $i\in I\setminus\{l\}$. Since $e_l=t\otimes E'_{-\theta}$
where $\theta$ is the highest root of $\dot \mg$, and
$\mg(A)_\b=t^{2k+1}\otimes \dot \mh_{\bar 1}$, then $[e_l,
\mg(A)_\b]=0$, and $d=\ad e_\b$.

If $A=A_{2l-1}^{(2)}$ and $\b=(2k+1)\d$,  we can  find
$e_\b\in\mg(A)_\b$ such that $\ad e_\b(e_i)=d(e_i)$ for all $i\in
I\setminus\{0,l\}$ because $\dim\mg(A)_\b=l-1$. Since $[e_l,
\mg(A)_\b]=[E'_l, \mg(A)_\b]=0$, we see that $(d-\ad
e_\b)(e_i)=\delta_{i,0}e_{\b+\a_0}$. If $e_{\b+\a_0}\ne0$, then
$e_{\b+\a_0}\in t^{2k+2}\otimes
F(E'_{-\theta^0}-E'_{-\bar\mu(\theta^0)})=t^{2k+2}\otimes \dot
\mg_{\bar 0}(-\theta_1)$ where $\theta_1$ is the highest short
root of $\dot \mg_{\bar 0}$ and $\theta^0$ is the root in $\dot
\D$ which induces the highest weight for $\dot \mg_{\bar 1}$ as a
$\dot \mg_{\bar 0}$-module, and
$E'_{-\theta^0}-E'_{-\bar\mu(\theta^0)}$ is the lowest weight
vector of $\dot \mg_{\bar 1}$. (see Section 8.3 at \cite{Kac}). We
use $i_0$ defined in the above discussion for the simple Lie
algebra $\dot \mg_{\bar 0}$ (of type $C_l$) defined in Theorem 8.3
at \cite{Kac}, and let $r_0=\a_0(\a_{i_0}^\vee)$. By computation
we see that $(\ad e_{i_0})^{r_0+1}e_0=0$ but $d((\ad
e_{i_0})^{r_0+1}e_0)=(\ad e_{i_0})^{r_0+1}d(e_0)=(\ad
e_{i_0})^{r_0+1}(e_{\b+\a_0})\ne0$ since $[f_{i_0},
(E'_{-\theta^0}-E'_{-\bar\mu(\theta^0)})]\ne0$. This is a
contradiction. Thus $e_{\b+\a_0}=0$, and $d=\ad e_\b$.

If $A=D_{l+1}^{(2)}$ and $\b=(2k+1)\d$,  we know that
$\dim\mg(A)_\b=1$, and $[e_i,\mg(A)_\b]=0$ for $i=1, 2, ..., l-1$.
We can  find $e_\b\in\mg(A)_\b$ such that $\ad e_\b(e_i)=d(e_i)$
for all $i\in I\setminus\{0\}$, i.e., $(d-\ad
e_\b)(e_i)=\delta_{i,0}e_{\b+\a_0}$. If $e_{\b+\a_0}\ne0$, then
$e_{\b+\a_0}\in t^{2k+2}\otimes
F(E'_{-\theta^0}-E'_{-\bar\mu(\theta^0)})=t^{2k+2}\otimes
\mg_{\bar 0}(-\theta_1)$ where $\theta_1$ is the highest short
root of $\dot \mg_{\bar 0}$ and $\theta^0$ is the root in $\dot
\D$ which induces the highest weight for $\dot \mg_{\bar 1}$ as a
$\dot \mg_{\bar 0}$-module, and
$E'_{-\theta^0}-E'_{-\bar\mu(\theta^0)}$ is the lowest weight
vector of $\dot \mg_{\bar 1}$. We use $i_0$ defined in the above
discussion for the simple Lie algebra $\dot \mg_{\bar 0}$ (of type
$B_l$) defined in Theorem 8.3 at \cite{Kac}, and let
$r_0=\a_0(\a_{i_0}^\vee)$. By computation we see that $(\ad
e_{i_0})^{r_0+1}e_0=0$ but $d((\ad e_{i_0})^{r_0+1}e_0)=(\ad
e_{i_0})^{r_0+1}d(e_0)= (\ad e_{i_0})^{r_0+1}(e_{\b+\a_0})\ne0$
since $[f_{i_0}, (E'_{-\theta^0}-E'_{-\bar\mu(\theta^0)})]\ne0$.
This is a contradiction. Thus $e_{\b+\a_0}=0$, and $d=\ad e_\b$.

If $A=E_{6}^{(2)}$ and $\b=(2k+1)\d$, we know that
$\dim\mg(A)_\b=2$, and $[e_i,\mg(A)_\b]=0$ for $i=3,4$. We can
choose $e_\b\in\mg(A)_\b$ such that $\ad e_\b(e_i)=d(e_i)$ for all
$i\in I\setminus\{0\}$, i.e., $(d-\ad
e_\b)(e_i)=\delta_{i,0}e_{\b+\a_0}$. If $e_{\b+\a_0}\ne0$, then
$e_{\b+\a_0}\in t^{2k+2}\otimes
F(E'_{-\theta^0}-E'_{-\bar\mu(\theta^0)})=t^{2k+2}\otimes \dot
\mg_{\bar 0}(-\theta_1)$ where $\theta_1$ is the highest short
root of $\dot \mg_{\bar 0}$ and $\theta^0$ is the root in $\dot
\D$ which induces the highest weight for $\dot \mg_{\bar 1}$ as a
$\dot \mg_{\bar 0}$-module, and
$E'_{-\theta^0}-E'_{-\bar\mu(\theta^0)}$ is the lowest weight
vector of $\dot \mg_{\bar 1}$. We use $i_0$ defined in the above
discussion for the simple Lie algebra $\dot \mg_{\bar 0}$ (of type
$F_4$) defined in Theorem 8.3 at \cite{Kac}, and let
$r_0=\a_0(\a_{i_0}^\vee)$. By computation we see that $(\ad
e_{i_0})^{r_0+1}e_0=0$ but $d((\ad e_{i_0})^{r_0+1}e_0)=(\ad
e_{i_0})^{r_0+1}d(e_0)\ne0$ since $[f_{i_0},
(E'_{-\theta^0}-E'_{-\bar\mu(\theta^0)})]\ne0$. This is a
contradiction. Thus $e_{\b+\a_0}=0$, and $d=\ad e_\b $.

Now assume that $A=D_{4}^{(3)}$,  and $\b=(3k+s)\d$ where $s=1$ or
$2$. Note that we used the primitive cubic root, $\eta$, of $1$,
in the construction of the affine Lie algebra, but this element is
not necessarily in $F$.
We know that $\dim\mg(A)_\b=1$, and
$[e_2,\mg(A)_\b]=0$. We can choose $e_\b\in\mg(A)_\b$ such that
$\ad e_\b(e_i)=d(e_i)$ for all $i\in I\setminus\{0\}$, i.e.,
$(d-\ad e_\b)(e_i)=\delta_{i,0}e_{\b+\a_0}$ where we can choose
$e_\b\in\mg(A)_\b$ such that $e_{\b+\a_0}=[e_{\b},e_{0}]$. Now we
replace $d$ by $d-\ad e_\b$.

If $e_{\b+\a_0}\ne0$ and $s=2$, then $e_{\b+\a_0}\in
t^{3k+3}\otimes F(E'_{-\theta^0}+E'_{-\bar\mu(\theta^0)}+
E'_{-\bar\mu^2(\theta^0)})=t^{3k+3}\otimes \dot \mg_{\bar
0}(-\theta_1)$ where $\theta_1$ is the highest short root of $\dot
\mg_{\bar 0}$, $\theta^0$ is the root in $\dot \D$ which induces
the highest weight for $\dot \mg_{\bar 1}$ as a $\dot \mg_{\bar
0}$-module,  $E'_{-\theta^0}+\eta ^2E'_{-\bar\mu(\theta^0)}+\eta
E'_{-\bar\mu^2(\theta^0)}$ is the lowest weight vector of $\dot
\mg_{\bar 1}$. We use $i_0$ defined in the above
discussion for the simple Lie algebra $\dot \mg_{\bar 0}$ (of type
$G_2$) defined in Theorem 8.3 at \cite{Kac}, and let
$r_0=\a_0(\a_{i_0}^\vee)$. By computation we see that $(\ad
e_{i_0})^{r_0+1}e_0=0$ but $d((\ad e_{i_0})^{r_0+1}e_0)=(\ad
e_{i_0})^{r_0+1}d(e_0)\ne0$ since
$[f_{i_0},e_{\b+\a_0}]=t^{3k+3}\otimes[f_{i_0},
E'_{-\theta^0}+E'_{-\bar\mu(\theta^0)}+
E'_{-\bar\mu^2(\theta^0)}]\ne0$. This is a contradiction. Thus
$e_{\b+\a_0}=0$, and $d=\ad e_\b$.

If $e_{\b+\a_0}\ne0$ and $s=1$, from $d([e_0,[e_0,e_1]])=0$ we
deduce that $[e_0,[[e_\b,e_0],e_1]]+[[e_\b,e_0],[e_0,e_1]]=0$ and
further $[e_0,[e_0,[e_\b,e_1]]]=0$. By concrete computation we
obtain that
$[e_0,[e_0,[e_\b,e_1]]]=t^{3k+3}\otimes[E'_{-\theta^0}+\eta
^2E'_{-\bar\mu(\theta^0)}+\eta E'_{-\bar\mu^2(\theta^0)},
[E'_{-\theta^0}+\eta ^2E'_{-\bar\mu(\theta^0)}+\eta
E'_{-\bar\mu^2(\theta^0)},E'_1+\eta^2E'_3+\eta E'_4]] \ne0$, since
$[e_\b,e_1] = t^{3k+1} \otimes (E'_1+\eta^2E'_3+\eta E'_4)$. This
is a contradiction. Thus $e_{\b+\a_0}=0$, and $d=\ad e_\b$.

Thus we exhausted all affine Lie algebras. This completes the
proof.\qed

Now we summarize our results in this section into the following.

\begin{thm} Suppose that $A$ is an indecomposable $(l+1) \times (l+1)$ GCM with $l \geq 1$.

(i) If $A$ is of finite type, then $\Der(\mn^+)$ is spanned by
{\rm $\ad(\mb^+)|_{\mn^+}$} and the following $l+1$ outer
derivations $d_i\ (i \in I)$ given by
$$d_i: e_{j}\mapsto
\delta_{i,j}e_{s_i(\theta)}\,\, (\in \mg_{s_i(\theta)})
\,\,\forall\,\, j\in I\mbox{;}$$

(ii) If $A$ is of affine type $X_{N}^{(r)}$, then $\Der(\mn^+)$ is
spanned by {\rm $\ad(\mb^+)|_{\mn^+}$} and the following outer
derivations $d_k\ (k \in \N)$ given by
$$d_k:
e_{i}\to \delta_{i,\e}e_{kr\delta+\a_\e}\,\, \,\,\forall\,\,
i\in I$$ where $e_{kr\delta+\a_\e}\in\mg_{kr\delta+\a_\e}$;

(iii) If $A$ is of indefinite type, then {\rm
$\Der(\mn^+)=\ad(\mb^+)|_{\mn^+}$} and {\rm
$\Der(\tn^+)=\ad(\tb^+)|_{\tn^+}$}.
\end{thm}

\noindent {\bf Remark 2}.\ When $A$ is of finite type,
$\Der(\mn^+)$ was   determined in (cf.~\cite{Ko}, \cite{LL2}).
When $A$ is of affine type, $\Der(\mn^+)$ was obtained in
(cf.~\cite{F1}, \cite{F2}).

\section{\textbf{ Automorphisms of $\mn^+$}}

In this section we assume that $A$ is an $(l+1)\times(l+1)$
indecomposable symmetrizable GCM. It is easy to see that the center
$Z(\mn^+) = 0$ if $\dim \mn^+ = \infty$ (this is true also for
nonsymmetrizable $A$). Denote $\bar\mb^+=\mb^+/\mc$. The main result
in this section is the following:

\begin{thm} Let $A$ be an indecomposable symmetrizable GCM of infinite type (that is, of affine type or
of indefinite type). Then,
$$\Aut(\mn^+) \simeq \Aut(\bar \mb^+) = (\Aut(A) \ltimes T) \ltimes U_+
= \Aut(A) \ltimes B_+\ .$$
\end{thm}

\noindent {\it Proof.} {\bf Case 1}:  $A$ is of indefinite type.

From Theorem 3.4(iii) we know that $\Der(\mn^+) =
\ad(\mb^+)|_{\mn^+}$. Hence, we have
$$\mn^+ \simeq \ad(\mn^+) \subset \Der(\mn^+) \simeq \bar \mb^+ = \mb^+/\mc\ .$$
If $\tau \in \Aut(\mn^+)$, then
$$\hat \tau : d \mapsto \tau \circ d \circ \tau^{-1}$$
gives $\hat \tau \in \Aut(\Der(\mn^+))$. Hence, we see that the
composition of maps (of the extension $\hat{\ }$ and the restriction
${\ }|_{\mn^+}$) defined by
$$\Aut(\mn^+) \overset{\widehat{\ }}{\longrightarrow} \Aut(\Der(\mn^+))
\overset{{\ }|_{\mn^+}}{\longrightarrow} \Aut(\mn^+)$$ is the
identity map: $\Aut(\mn^+) \longrightarrow \Aut(\mn^+)$.

We consider $\bar \mh = \mh/\mc \subset \bar \mg = \mg/\mc$ with
$\mn^+ \subset \bar \mb^+ = \bar \mh \oplus \mn^+ \subset \bar
\mg$. Let $\tau \in \Aut(\mn^+)$, and extend $\tau$ to $\hat \tau
\in \Aut(\Der(\mn^+)) = \Aut(\bar \mb^+)$. Put $\bar {\mathfrak t}
= \hat \tau (\mh) \subset \bar \mb^+ \subset \bar \mg$. Using the
same method as in Section 2, we find that $\bar {\mathfrak t}$ is
a split Cartan subalgebra of $\bar \mg$. By the conjugacy result for
split Cartan subalgebras of $\bar \mg$ (cf.~\cite{PK}), we will
determine $\Aut(\bar \mb^+)$ as follows. Again using the method as
in Section 2, there is an element $\bar u \in \bar U_+ =
U_+|_{\bar \mn^+} \simeq U_+$ such that $\bar u(\bar {\mathfrak
t}) = \bar \mh$. Replacing $\tau$ by $\bar u\tau$, we can assume
$\hat \tau(\mh) = \mh$. Taking modulo $\Aut(A)$ as in Section 2,
we can also assume $\hat \tau(\bar \mg_{\a_i}) = \bar \mg_{\a_i}$
for all $i \in I$. Then, modulo $T$, we can assume $\hat \tau(\bar
e_i) = \bar e_i$ for all $i \in I$. Therefore, we can reach $\hat
\tau(\bar \a_i^*) = \bar \a_i^*$ for all $i \in I$, where $\bar
\a_i^* = \a_i^*\ (\mbox{modulo}\ \mc)$ in $\bar \mh$. Thus, we
have $\Aut(\bar \mb^+) = (\Aut(A) \ltimes T) \ltimes U_+ = \Aut(A) \ltimes B_+$. Then, we
can easily confirm that $\Aut(\bar \mb^+)$ is isomorphic to
$\Aut(\mn^+) = \Aut(\bar \mb^+)|_{\mn^+}$.

{\bf Case 2}:  $A$ is of affine type.

Let $D = \Der(\mn^+)$. Then, $D = D_0 \oplus (\oplus_{\a \in
\D_+}\, D_\a)$ by Section 3. Using a standard loop algebra
realization and Section 4 (also see \cite{F2}), we have
$$D = \mathcal{V} \oplus \ad(\mh' \oplus \mn^+)|_{\mn^+} = \mathcal{V}_+ \oplus
\ad(\mb^+)|_{\mn^+}\ ,$$
where $\mathcal{V}$ and $\mathcal{V}_+$ are defined by
$$\mathcal{V} = \oplus_{p \geq 0,\ p \equiv 0\ ({\rm mod}\, r)}\, F\partial_p\ ,\quad
\mathcal{V}_+ = \oplus_{p > 0,\ p \equiv 0\, ({\rm mod}\, r)}\, F\partial_p\ ,\quad
\partial_p = t^p\cdot t\frac{d}{dt}$$
with Virasoro type relations (cf.~\cite{Wa})
$$[\partial_p,\partial_q] = (q-p)\partial_{p+q}$$
for all possible $p,q$.
We note that $\ad(\mb^+)|_{\mn^+} \simeq \bar \mb^+$ is a subalgebra of $D$, and that
$\ad(\mh' \oplus \mn^+)|_{\mn^+}$ is an ideal of $D$.
Let
$$\Theta = \{ d \in D \,|\, \ad(d)\ \mbox{is locally finite on}\ D \}\ .$$
Since all of $D_0 = \ad(\mh)|_{\mn^+}$ and
$D_\a = \ad(\mg_\a)|_{\mn^+}$ for $\a \in \D_+^{\rm re}$ are contained in $\Theta$,
we see that $\ad(\mb^+)|_{\mn^+}$ is contained in the subalgebra,
$\langle \Theta \rangle$, of $D$ generated by $\Theta$.
Suppose $\ad(\mb^+)|_{\mn^+} \ne \langle \Theta \rangle$.
Then, we can find $d' \in \Theta$ such that $d' \not\in \ad(\mb^+)|_{\mn^+}$.
Put $\bar D = D/(\ad(\mh' \oplus \mn^+)|_{\mn^+}) \simeq \mathcal{V}$.
Using Virasoro type relations, for a nonzero element
$\partial = \sum_{i=1}^{k}\, a_{p_i}\partial_{p_i} \in \mathcal{V}$ with
$$0 \leq p_1 < p_2 < \cdots < p_k,\quad p_i \equiv 0\, ({\rm mod}\ r),\quad
a_{p_k} \ne 0,\quad k \geq 1\ ,$$
we easily see that
$\ad(\partial)$ is locally finite on $\mathcal{V}$ if and only if $k = 1$ and $p_k = 0$.
Hence, $\bar d'\, (\in \bar D)$ must be in $F\partial_0 \subset \mathcal{V}$, which means
$d' \in \ad(\mb^+)|_{\mn^+}$.
This is a contradiction. Therefore, we obtain
$\ad(\mb^+)|_{\mn^+} = \langle \Theta \rangle$.
In particular, $\ad(\mb^+)|_{\mn^+}$ is a characteristic subalgebra of $D$.

Now, we can return to the same setting as in case of indefinite
type. Namely, let $\tau \in \Aut(\mn^+)$. Then, we can extend
$\tau$ to an automorphism $\hat \tau$ of $D = \Der(\mn^+)$ as
before. Note that $\ad(\mb^+)|_{\mn^+}$ is invariant under the
action of $\Aut(D)$, since $\ad(\mb^+)|_{\mn^+}$ is a
characteristic subalgebra of $D$. Hence, $\hat \tau$ induces an
automorphism of $\ad(\mb^+)|_{\mn^+}$, again called $\hat \tau$.
Using the same method as in case of indefinite type, we have
$\Aut(\ad(\mb^+)|_{\mn^+}) \simeq \Aut(\bar \mb^+) = (\Aut(A)
\ltimes T) \ltimes U_+ = \Aut(A) \ltimes B_+$. Then, taking a
restriction of $\hat \tau$ to $\ad(\mn^+) \simeq \mn^+$, we can
obtain the result in Theorem 4.1.\qed

To make the result in Theorem 4.1 complete, we need to give
$\Aut(\mn^+)$ when $A$ is of finite type but not of $A_1$ (since
it is trivial).

From Theorem 3.1(i), we know that $D=\Der(\mn^+)$ is spanned by {\rm
{\rm $\ad(\mb^+)|_{\mn^+}$} and the $l+1$ outer derivations $d_i\ (i
\in I)$.
Set
$$D_0 = \ad(\mh)|_{\mn^+}\ ,
\,\,\,
D' = \oplus_{i \in I}\, Fd_i\ ,
\,\,\,
D'' = \ad(\mn^+) \simeq \mn^+/\mg_\theta\ ,
\,\,\,
D_+ = D' \oplus D''\ .
$$
Then, $D = \oplus_{\b \in Q}\, D_\b = D_0 \oplus D_+$,
and $D''$ is an ideal of $D$.
We note that
$[\ad h,d_\b] = \b(h)d_\b$ for all $h \in \mh$ and $d_\b \in D_\b$, and that
$[d,\ad x] = \ad d(x)$ for all $d \in D$ and $x \in \mn^+$.

From Theorem 3.4, we see that if $A
\ne {\rm A}_2, {\rm B}_2$, then $D_+ = \{ d \in D\, |\, d(\mn^+)
\subset [\mn^+,\mn^+] \}$. In fact, ${\rm A}_2$ is special, so we need to
handle it separately. But ${\rm B}_2$ is not necessarily excluded,
since our proof for genral case also works in case of type ${\rm B}_2$.

Let $\exp(\ad D'')$ be the subgroup of $\Aut(D)$
generated by $\exp(\ad x)$ for all $x\in D''$. Then $\exp D''
\simeq U_+/Z(U_+)$ as automorphisms of $\mn^+$, where $Z(U_+)$ is
the center of $U_+$. An automorphism $\tau$ of $\mn^+$ induces an
automorphism, called $\bar \tau$, of $\ad(\mn^+) \simeq
\mn^+/\mg_\theta$, which is denoted by $\zeta:\Aut(\mn^+)
\longrightarrow \Aut(\mn^+/\mg_\theta)$, $\tau\to \bar\tau$.
We note that if $A \ne {\rm A}_2$ then
$[D,D]=[D,D_+]=D_+$ and $[D',D']=0$, but
we see that if $A = {\rm A}_2$ then
$[d_0,d_1] = \ad(-\a_0^*+\a_1^*) \in D_0$.

Now we can obtain the following:

\begin{thm} Let $A$ be an $(l+1)\times(l+1)$ indecomposable GCM of finite type
with $l \geq 1$. Then we have
$$\Aut(\mn^+)
= \langle \Aut(A),\ T,\ \exp(Fd_i)\, |\, i \in I \rangle \ltimes \exp D''\ ,$$
more precisely,
$\Aut(\mn^+)
= (\Aut(A) \ltimes T) \ltimes \exp D_+ \simeq (\Aut(A)
\ltimes (T \ltimes \exp D')) \ltimes U_+/Z(U_+)$ if $A \ne {\rm
A}_2$,
while $\Aut(\mn^+)\simeq {\rm GL}(2) \ltimes F^2$ if $A =
{\rm A}_2$.
\end{thm}

\noindent {\it Proof.}  In case of $A = {\rm A}_2$, we may assume
$d_0(e_j) = \d_{0,j}e_1$ and $d_1(e_j) = \d_{1,j}e_0$ for $j \in I
= \{ 0,1 \}$. Since $d_0, d_1, [d_0, d_1]$ form a copy of
$\mathfrak{s}\mathfrak{l}_2$, we see that $D$ is not solvable. In
this case, we can directly compute $\Aut(\mn^+)$ as follows. Put
$x = e_0,y = e_1, z = [x,y] = [e_0,e_1]$, as a basis of $\mn^+$ (a
$3$ dimensional Heisenberg Lie algebra). Let $\tau \in
\Aut(\mn^+)$, and write $\tau(x) = a_1x+a_2y+a_3z$ and
$\tau(y)=b_1x+b_2y+b_3z$. Then, $a_1b_2-b_1a_2 \ne 0$ and $\tau(z)
= (a_1b_2-b_1a_2)z$. Conversely this condition gives an
automorphism of $\mn^+$. With respect to a basis $x,y,z$ of
$\mn^+$, the automorphism group $\Aut(\mn^+)$ can be realized as
$$\left\{ \left.
\left(
\begin{array}{ccc}
a_1 & b_1 & 0\\ a_2 & b_2 & 0\\ a_3 & b_3 & a_1b_2-b_1a_2
\end{array}
\right)\, \right|\, a_1,a_2,a_3,b_1,b_2,b_3 \in F,\, a_1b_2-b_1a_2
\ne 0\, \right\}\ \le {\rm GL}(3).$$  Furthermore, we have
$$\exp D'' \leftrightarrow \left\{ \left. \left(
\begin{array}{ccc} 1&0&0\\ 0&1&0\\ a&b&1
\end{array}\right)\, \right|\, a,b\in F\right\}
\simeq F^2,$$ $$T\leftrightarrow \left\{ \left. \left(
\begin{array}{ccc} a&0&0\\ 0&b&0\\ 0&0&ab
\end{array}\right)\, \right|\, a,b \in F^\times \right\}
\simeq {F^\times}^2,$$
$$\exp(Fd_0) \leftrightarrow \left\{ \left. \left(
\begin{array}{ccc} 1&0&0\\ a_2&1&0\\ 0&0&1
\end{array}\right)\, \right|\, a_2\in F\right\}
\simeq F,$$
$$\exp(Fd_1) \leftrightarrow \left\{ \left. \left(
\begin{array}{ccc} 1&b_1&0\\ 0&1&0\\ 0&0&1
\end{array}\right)\, \right|\, b_1\in F\right\}
\simeq F,$$
$$\left( \exp(D' \oplus F[d_0,d_1]) \leftrightarrow
\left\{ \left. \left(
\begin{array}{ccc} a_1&b_1&0\\ a_2&b_2&0\\ 0&0&1
\end{array}\right)\, \right|\,
\begin{array}{l}
a_1,a_2,b_1,b_2 \in F\\
a_1b_2-a_2b_1=1
\end{array}
\right\} \simeq SL(2)\ {\rm if}\ F=\C, \right)
$$
$$\Aut(A) \leftrightarrow \left\langle \left(
\begin{array}{ccc} 0&1&0\\ 1&0&0\\ 0&0&-1
\end{array}\right) \right\rangle \simeq \Z/2\Z.$$

Hence, the theorem follows in this case. We note that if $F=\Bbb C$ then
$\exp(D_0 \oplus D') \simeq {\rm GL}(2)$ (cf.~Remark 3, later).

Suppose $A \ne {\rm A}_2$. We note that $D$ is solvable and $D_+ =
[D,D]$, which is the nilradical of $D$.
In this case,
$D_+$ is an ideal of $D$, and $D'$ is
an abelian subalgebra of $D$.
We also define $\exp(\ad D')$ and $\exp(\ad D_+)$ by
$\exp(\ad D') = \langle \exp(\ad x')\, |\, x' \in D' \rangle$ and
$\exp(\ad D_+) = \langle \exp(\ad x_+)\, |\, x_+ \in D_+ \rangle$
respectively as subgroups of $\Aut(D)$.
Using the following well-known formulas
$$\varphi\exp(\ad d)\varphi^{-1}=\exp(\ad\varphi(d)),\,\,\forall \,\,d\in D, \varphi\in \Aut(D),\,\, {\rm and}$$
$$\exp(\ad d)\exp(\ad d')\exp(-\ad d)=\exp(\ad[d,d']),\,\,\forall \,\,d, d'\in D,$$ we can easily verify that  $\exp(\ad D')$ is an abelian subgroup
of $\Aut(D)$ with $\exp(\ad D') \cong F^{l+1}$, and $\exp(\ad D_+)$ is a
normal subgroup of $\Aut(D)$, also $\exp(\ad D'')$  is a normal subgroup
of $\exp(\ad D_+)$, and further
$$\exp(\ad D_+) = \exp(\ad D') \ltimes \exp(\ad D'') \subset
\Aut(D).$$

Let $\tau \in \Aut(\mn^+)$. We define $\hat \tau \in \Aut(D)$ by
$\hat \tau(d) = \tau \circ d \circ \tau^{-1}$. From $\hat\tau(\ad
x)=\ad(\tau(x))$, clearly, $\hat \tau (D'')=D''$. Put
$\mathfrak{t} = \hat \tau(D_0)$. Then, both $D_0$ and
$\mathfrak{t}$ are split Cartan subalgebras of $D$. Since $D_+$ is
the nilpotent radical of $D$, we can find $\hat \tau' \in \exp(\ad
D_+)$ such that $\hat \tau'(\mathfrak{t}) = D_0$ (see Ex.7.6,
\cite{MP}). Put $\hat \tau'' = \hat \tau' \hat \tau \in \Aut(D)$.
Then, $\hat \tau''(D_0) = D_0$ and $\hat \tau''(D'') = \hat
\tau'(\hat \tau(D'')) = \hat \tau'(D'') = D''$. Hence, $\hat
\tau''$ induces an automorphism of $D_0 \oplus D'' =
\ad(\mb^+)|_{\mn^+} \simeq \mb^+/\mg_\theta = \mh \oplus
(\mn^+/\mg_\theta)$. Using a standard argument as in Section 2, we
see that $\hat \tau''(\a_i^\vee) = \a_i^\vee$ and $\hat
\tau''(\bar e_i) = \bar e_i$ for all $i \in I$ modulo the action
of $\Aut(A) \ltimes T$, where $\bar e_i = e_i \mod \mg_\theta$.
That is, $(\sigma \xi \hat \tau' \hat \tau)|_{D_0 \oplus D''} = 1$
for some $\sigma \in \Aut(A)$ and $\xi \in T$. In particular, we
have $(\sigma \xi \hat \tau' \hat \tau)|_{D''} = 1$. Hence,
$\sigma|_{D''} \circ \xi|_{D''} \circ \hat \tau'|_{D''} \circ \hat
\tau|_{D''} = 1$. Since $\hat \tau|_{D''} = \bar \tau$ and
$\Aut(A)\, T\, \exp D_+ \subset \Aut(\mn^+)$, we obtain
${\rm Im}(\zeta) = \Aut(A)\, T\, \exp(\ad D_+)|_{D''}$. On
the other hand, we have ${\rm Ker}(\zeta) = \langle \tau_i(z)\,
|\, i \in I,\ z \in \mg_\theta \rangle \simeq F^{l+1}$, where $z$
can be any element of $\mg_\theta$ and $\tau_i(z) \in \Aut(\mn^+)$
is given by $\tau_i(z):e_j \mapsto e_j+\d_{i,j}z\ (j \in I)$. If
$\theta(\a_i^\vee) > 0$, then $\tau_i(\mg_\theta) = \exp(\ad\
\mg_{\theta-\a_i})$. If $\theta(\a_i^\vee) = 0$, then
$\tau_i(\mg_\theta) = \exp(Fd_i)$. Therefore, we see that ${\rm
Ker}(\zeta) \subset \exp D_+$. Hence, we have $\Aut(\mn^+) =
(\Aut(A) \ltimes T) \ltimes \exp D_+ \simeq (\Aut(A) \ltimes (T
\ltimes \exp D')) \ltimes U_+/Z(U_+)$. \qed

{\bf Remark 3.} Note that if $F=\C$, then $T=\exp(\ad\
\mh)=\exp(D_0)$. Consequently, $\Aut(\mn^+) =\Aut(A)\ltimes \exp(D)$
if $A$ is of finite type but not $A_2$, and $ \Aut(\mn^+)= \exp D$
if $A=A_2$.

\section{\textbf{ Derivations of $\mb^+$ and $\tb^+$}}

In this section, a generalized Cartan matrix $A$ is not
necessarily symmetrizable. We will determine $\Der(\mb^+)$ and
$\Der(\tb^+)$ in an elementary way. The result in this section is:

\begin{thm}
$\Der(\mb^+) = \hom(\mh,\mc) \oplus {\rm ad}(\mb^+)$, and
$\Der(\tb^+) = \hom(\mh,\mc) \oplus {\rm ad}(\tb^+).$
\end{thm}

\noindent {\bf Proof.} Suppose $\bar F$ is an extension field of
$F$. For any Lie algebra $L$ over $F$, we denote $\Der_{\bar F}(\bar
F\otimes _F L)$ to be the derivation algebra of the Lie algebra
$\bar F\otimes _F L$. It is well-known that $\Der_{\bar F}(\bar
F\otimes L)=\bar F\otimes _F\Der(L)$ if $L$ is finitely generated.

Indeed, take a basis $\{1, a_j\, |\, j\in J\}$ for $\bar F/F$ and
let $d\in \Der_{\bar F}(\bar F\otimes _F L)$.  For any $x\in L$,
write $d(x)=x_0+\sum_{j\in J}a_jx_j$ where the sum is finite for
each $x$. For each $j\in J \cup \{ 0 \}$, define $d_j: L\to L,
x\mapsto x_j$. For any $x,y\in L$, from $$d_0([x,y])+\sum_{j\in
J}a_jd_j([x,y])=d([x,y])=[d(x),y]+[x,d(y)]$$
$$=([d_0(x),y]+[x,d_0(y)])+\sum_{j\in J}a_j([d_j(x),y]+[x,d_j(y)]),$$
we see that $d_j\in \Der(L)$.  Since $L$ is finitely generated, $J'
= \{ j \in J\, |\, d_j\ne0\
 \}$  is finite. Clearly, $d = d_0 + \sum_{j \in J'}\, a_j \otimes
d_j \in \bar F \otimes_F \Der(L)$. Then $\Der_{\bar F}(\bar F\otimes
L)=\bar F\otimes _F\Der(L)$.

Thus we may assume that $F$ is algebraically closed in the statement
of the theorem.

Since $F$ is algebraically closed, we can fix $h_0\in \mh$ so that
for any $\b_1,\b_2\in Q$, $$\b_1(h_0)=\b_2(h_0)\,\,{\rm iff} \,\,
\b_1=\b_2.$$

Let $d \in \Der(\mb^+)$ be a derivation of $\mb^+$. Here and
afterwards, $\mb^+$ denotes $\mb^+$ and $\tb^+$. We may assume
that
$$d(h_0)=h_1+\sum_{\a\in\Delta^+}x_\a,$$ where $h_1\in\mh,
x_\a\in\mb^+_\a$ and the sum on the right-hand side is finite.
Replacing $d$ by $d-\sum_{\a\in\Delta^+}\a(h_0)^{-1}\ad x_\a$, we
may assume that $d(h_0)=h_1\in\mh$. Since $0 = d([h_0,h]) =
[d(h_0),h] + [h_0,d(h)] = [h_0,d(h)]$ and $d(h) \in C_{\mb^+}(h_0)
= \mh$ for all $h \in \mh$ we deduce $d(\mh)\subset \mh$.

For any nonzero $x_\a\in\mb^+_\a$, we have
$$\a(h_0)d(x_\a)=d([h_0,x_\a])=[h_1,x_\a]+[h_0,d(x_\a)]=\a(h_1)x_\a+[h_0,d(x_\a)],$$
which yields $\a(h_1)=0$ and $d(x_\a)\in\mb^+_\a$ for all
$\a\in\Delta^+$. Thus $h_1\in \mc$ and $d(e_i)\in\mb^+_{\a_i}$ for
all $i\in I$. Thus we can find $h_2\in\mh$ such that $(d-\ad
h_2)(e_i)=0$ for all $i\in I$. Replacing $d$ by this $d-\ad h_2$,
we may assume that $d(\mn_+)=0$ and $d(\mh)\subset \mh$.

From $d([h,e_i])=0$ for all $h\in \mh$ and all $i\in I$, we deduce
that $d(h)\in\mc$ for all $h\in\mh$. Thus $d(\mh)\subset\mc$. This
completes the proof.\qed

{\bf Remark 4.} If $A$ is of finite type, the result in Theorem
5.1 was obtained in \cite{LL1} with a different approach, also see
\cite{HS}, \cite{Ko}, \cite{Fu} from Lie algebra cohomology view of point.

One can use the method established in the proof of Theorem 5.1 to show
the following theorem.

\begin{thm}
$\Der(\mg(A)) = \hom(\mh'',\mc) \oplus {\rm ad}(\mg(A))$ and
$\Der(L(A)) = \hom(\mh'',\mc) \oplus {\rm ad}(L(A))$.
\end{thm}

\noindent {\bf Proof.} Let $d \in Der(\mg)$ and take $h_0 \in \mh$
as in the proof of Theorem 5.1. We can assume that $d(h_0) = h_1 \in
\mh$ and $d(\mh) \subset \mh$ modulo $\ad(\mn^-) \oplus \ad(\mn^+)$.
Then, we can also reach $h_1 \in \mc$ and $d(\mg_{\pm \a_i}) \subset
\mg_{\pm \a_i}$ for all $i \in I$. For each $h \in \mh$ and $i \in
I$, we obtain that $[h,e_i] = \a_i(h)e_i$ implies $\a_i(h)d(e_i) =
d([h,e_i]) = [d(h),e_i] + [h,d(e_i)] = \a_i(d(h))e_i +
\a_i(h)d(e_i)$ and $\a_i(d(h))e_i = 0$. Hence, $\a_i(d(h))=0$ for
all $i \in I$, and $d(h) \in \mc$ for all $h \in \mh$. If $d(e_i) =
\lambda_ie_i$ and $d(f_i) = \mu_if_i$, then $d(\a_i^\vee) =
d([e_i,f_i]) = [d(e_i),f_i] + [e_i,d(f_i)] =
(\lambda_i+\mu_i)\a_i^\vee \in \mc$. Therefore, $\lambda_i + \mu_i =
0$. Then, there is $h_2 \in \mh$ such that $(d-\ad h_2)(e_i) =
(d-\ad h_2)(f_i) = 0$ for all $i \in I$. Hence, we can assume
$d(\mg') = 0$ modulo $\ad(\mg)$ for our purpose. Thus, we have
$\Der(\mg) \subset \hom(\mh'',\mc) \oplus {\rm ad}(\mg)$. Since the
other inclusion is trivial, we obtain the desired result for
$\Der(\mg(A))$.

The above argument is valid also for $L(A)$ with slight
modifications.\qed

{\bf Remark 5.} $\Der(\mg'/\mc)$ is determined by S. Berman for
non-affine case (cf.~\cite{B1}) and $\Der(\mg)$ is determined in
\cite{MP} (Sect.4, Prop.15, p.329) for finitely displayed
contragredient Lie algebras using different approach, while L.~J.
Santharoubane obtained $H^2(\mg,F) = \wedge^2 (\mh'')^*$ for affine
case (cf.~\cite{Sa}). According to Proposition 16.4 in \cite{N},
Theorems 5.1 and 5.2 can also be proved by computing degree $0$
derivations.

{\bf Remark 6.} The automorphism group $\mg'/\mc$ for symmetrizable
GCM $A$ was described by Peterson and Kac in \cite{PK}. We can find
that outer derivations $\hom(\mh'',\mc) \simeq \mc^{m'}$ are
corresponding to outer automorphisms $\Gamma_0(\mg) \simeq \mc^{m'}$
(cf.~Theorem 2.2, Theorem 5.2). For $\mb^+$, outer derivations
$\hom(\mh,\mc) \simeq {\rm End}(\mc) \oplus \mc^{l+1}$ are
corresponding to outer automorphisms $\Gamma_0(\mb^+) \simeq {\rm
GL}(\mc) \ltimes \mc^{l+1}$ (cf.~Theorem 2.1, Theorem 5.1).
Similarly we can also find, for $\mn^+$, a suitable correspondence
between locally nilpotent (or finite) outer derivations in
$\Der(\mn^+)$ and certain outer automorphisms in $\Aut(\mn^+)$ in
each case (cf.~Theorem 3.4, Theorem 4.1, Theorem 4.2). Finally we
note that $\dim H^1(\mn^+,\mn^+) = 2(l+1)$ if $A$ is of finite type
(assumed $A \ne {\rm A}_1$); $\dim H^1(\mn^+,\mn^+) = \infty$ if $A$
is of affine type; $\dim H^1(\mn^+,\mn^+) = l+1$ if $A$ is of
indefinite type.

\vskip 5pt {\bf Acknowledgement.} The authors would like to thank
Prof.\ Erhard Neher for providing \cite{N} and valuable suggestions
to improve the paper.

{\bf Addresses:} \

Institute of Mathematics, University of Tsukuba, Tsukuba, Ibaraki,
305-8571 Japan. Email: morita@math.tsukuba.ac.jp

Department of Mathematics, Wilfrid Laurier University, Waterloo,
ON, Canada N2L 3C5, and Academy of Mathematics and System
Sciences, Chinese Academy of Sciences, Beijing 100190, P. R.
China. Email: kzhao@wlu.ca

\enddocument